\documentclass{amsart}
\usepackage{amsmath}
\usepackage{amssymb}
\usepackage{amsthm}
\usepackage{amsfonts}

\def\ad{\mbox{ad} }
\def\h{\mbox{ht} }
\def\Ad{\mbox{Ad } }
\def\End{\mathop{\textup{End}} }
\def\Hom{\mathop{\textup{Hom}} }
\def\rk{\mathop{\textup{rk}} }
\def\sp{\mathop{\textup{span}}}
\def\M{{\mathcal M}}
\def\H{{\mathcal H}}

\def\rightarrows{\rightarrow\!\!\!\!\rightarrow}
\def\downarrows{\downarrow  \hspace{-.8em}\raisebox{-.25em}{$\downarrow$}}

\newcommand\p[1]{\rho_{#1}^*}
\newcommand\g[1]{\mathfrak{g}_{#1}}
\newcommand\n[1]{\mathfrak{n}_{#1}}

\newcommand\C[1]{{\mathcal C}_{#1}}
\newtheorem{theorem}{Theorem}
\newtheorem*{thmnonumber}{Theorem}
\newtheorem{lemma}[theorem]{Lemma}
\newtheorem{proposition}[theorem]{Proposition}
\newtheorem{definition}[theorem]{Definition}
\newtheorem{corollary}[theorem]{Corollary}

\begin{document}

%%%for article style 
\title[Decomposing Hessenberg varieties]{Decomposing 
Hessenberg varieties over Classical Groups}

%%%for puthesis style\title{Decomposing Hessenberg varieties over Classical Groups}

\author{Julianna S.~Tymoczko}
%%%for article style
\email{tymoczko@umich.edu}

%% *** Jared's comments: too many centered equations
%% ***                  repeated defns in sects 1 and 2!! 
%% ***                  examples section
%%%for article style \begin{document}
%%%for article style 
\begin{abstract}
%%%for puthesis style \abstract{
Hessenberg 
varieties are a family of subvarieties of the flag variety, 
including
the Springer fibers, the Peterson variety, and the entire flag variety itself.
  The
seminal example arises from a problem in numerical analysis and consists for a fixed linear
operator $M$ of the full flags $V_1 \subsetneq V_2 \ldots  \subsetneq V_n$ in 
$GL_n$
with $M V_i \subseteq V_{i+1}$ for all $i$. 

In this paper I show that all Hessenberg varieties in type $A_n$ and 
semisimple and regular
nilpotent Hessenberg varieties in types $B_n$,$C_n$, and $D_n$ 
can be paved by affine spaces.  
Moreover,
this paving is the intersection of 
a particular Bruhat decomposition with the Hessenberg 
variety.  In type $A_n$, an equivalent description of the cells of
the paving in terms of certain fillings of a Young diagram
can be used to compute the Betti numbers of Hessenberg varieties.
As an example, I show that the Poincare polynomial of the Peterson variety in $A_n$
is $\sum_{i =0}^{n-1} \binom{n-1}{i} x^{2i}.$
%%%for puthesis style}
%%%for article style
\end{abstract}

\maketitle

%%%for article style
\thanks{{\bf Acknowledgements.}
%%%for puthesis style\acknowledgements{
I would like to thank Gil Kalai, David Kazhdan, Arun 
Ram, and Eric Sommers
for many helpful conversations and suggestions.  Discussions at different 
stages with Jared
Anderson, Andrew Booker, Emma Carberry, Henry Cohn, Ketan Delal, 
Jordan Ellenberg, Miranda Hodgson, Allison Klein, and David Nadler
provided me with useful suggestions, insights, and morale.  I 
particularly want to thank Jared Anderson, Henry Cohn, and Eric Sommers for
reading early drafts of this thesis and offering corrections and
improvements.

I am grateful for the support and entertainment given me 
continually by the extended Tymoczko clan:
Maria Tymoczko, 
Marlene Wong,
Alexei Tymoczko,
Dmitri Tymoczko,
Elisabeth Camp,
Molly Donohue, and
Misha Kazhdan.  I am grateful also to Marshall Poe.

Robert MacPherson taught me how to think about mathematics.  This thesis
would have been impossible without him.}
%%%\end{acknowledgements}
%%%for puthesis style}
%%%\dedication{**FILL ME IN!**}

\section{Introduction}

Hessenberg varieties form a large class of subvarieties of the flag variety,
many examples of which have been of great importance to geometers, 
representation theorists, combinatorists, and numerical analysts, among others.  
In this paper I describe the basic topology of many Hessenberg
varieties.

Given a Lie algebra $\g{}$ with a Borel subalgebra $\mathfrak{b}$, a 
Hessenberg space $H$ is a $\mathfrak{b}$-submodule of $\g{}$ which contains
$\mathfrak{b}$.  For a fixed element $M$ in $\g{}$, we can consider
the elements $g$ in an associated linear algebraic group $G$ such 
that $\Ad g^{-1}(M)$ lies in $H$.
This gives a subset $G(M,H)$ of the linear algebraic group.  
Since $H$ is closed under conjugation by the elements of the Borel subgroup
$B$ which corresponds to $\mathfrak{b}$,
the subset $G(M,H)$ is closed under
right multiplication by elements of $B$.
Thus the image of $G(M,H)$ in the flag variety $G/B$ is
a closed subvariety $\H(M,H)$ of $G/B$.  This subvariety
$\H(M,H)$ is the Hessenberg variety of $M$ and $H$.  

Hessenberg varieties as such were introduced by De Mari, Procesi, and
Shayman in \cite{MPS}.
De Mari and Shayman were first motivated to study these
spaces because of a question in numerical analysis related to efficient
computation of the eigenvalues and eigenspaces of the operator $M$.
Given certain $H$, the space $\H(M,H)$ parametrizes the bases with respect to which
the operator $M$ can be efficiently diagonalized via the QR-algorithm 
\cite{MS}.
In \cite{MPS}, the authors provided a cell decomposition of 
$\H(M,H)$ when $M$ is regular
semisimple by using a natural torus action that exists for those $M$.  They 
observed
that when $H$ is generated by $\mathfrak{b}$ as well as the simple negative root spaces
then $\H(M,H)$ is the toric variety associated to the decomposition into Weyl
chambers.  This space is combinatorially interesting as well, since
its Betti numbers generalize Eulerian descents of a permutation and can be 
used to give generating functions for several permutation statistics \cite{F}.

For entirely different reasons, several major examples of 
nilpotent Hessenberg varieties have been intensely studied recently.  
Springer 
initiated this research when he discovered an amazing connection 
between the cohomology of the Springer fibers and the irreducible 
representations of the Weyl group \cite{S}.  These Springer fibers are
in fact the nilpotent
Hessenberg varieties $\H(N,\mathfrak{b})$.
Springer's original proof was algebraic but later work 
expanded on the geometric nature of the results, including \cite{BM},
\cite{CG}, \cite{Ho}, \cite{KL}, and \cite{L}, among others.
Spaltenstein identified the components of type-$A_n$ Springer fibers 
$\H(N,\mathfrak{b})$ and proved they were equidimensional
and then extended the proof of equidimensionality to general Springer fibers 
(see \cite{Sp1}, \cite{Sp2}).  Shimomura partitioned type-$A_n$ Springer fibers 
into 
affine spaces in a manner similar to that used here \cite{S1}. 
Spaltenstein further 
showed that there is a Schubert decomposition whose intersection with
the Springer fibers gives a paving by 
affines for type $A_n$ in \cite[section II.5]{Sp3} 
and Shimomura 
extended this to apply to partial flag varieties \cite{S2}.
%%Each author generalized
%%his results in \cite{Sp2} and \cite{S2} respectively.
%% *** check restrictions: Sp2 for any unipotents?, Sp1 and S1 for GL,
%% *** S2 for which groups?
Spaltenstein also gave a combinatorial description of
the cells used in this paving \cite{Sp1}.  Several years later,
De Concini, Lusztig, and Procesi
provided a paving by affines of the Springer fibers for all classical 
types by reducing to the case of Springer fibers of 
distinguished nilpotents \cite{CLP}.  In their work,
Borho and MacPherson generalized Springer fibers to the larger class of
nilpotent Hessenberg varieties given by $\H(N,\mathfrak{p})$ for each
parabolic subalgebra $\mathfrak{p}$.  They 
showed that the intersection cohomology
of these Hessenberg varieties also could be viewed as 
representations of the Weyl group \cite{BM}.

More recently still, Peterson defined the Peterson variety, which
 plays a role in quantum cohomology and whose totally
positive part has interesting properties.  
The Peterson variety is the nilpotent Hessenberg variety $\H(N,H)$
when $N$ is regular and $H$ is the Hessenberg space generated by $\mathfrak{b}$
together with the simple negative root spaces.  Kostant showed that the 
coordinate ring of a particular open affine subvariety of
the Peterson variety coincides with the quantum cohomology of the flag variety
\cite{Ko}.
Rietsch has shown that the totally nonnegative part of the Peterson 
variety $\H(N,H)$ is homeomorphic to the totally nonnegative part of
Givental's critical point locus for the mirror symmetric family for the
flag variety
\cite{R}.
Research into the Peterson variety is ongoing.

In the rest of this paper I describe general Hessenberg varieties
and then give a paving by affine spaces for all Hessenberg varieties in
type $A_n$ as well as semisimple and 
regular nilpotent Hessenberg varieties in the other
classical types.  The main theorem for regular nilpotent Hessenberg varieties
is similar to that for semisimple Hessenberg varieties in classical types
and general Hessenberg varieties in type $A_n$.

\begin{thmnonumber}
Fix a regular nilpotent element $N$, let $\mathfrak{b}$ be the unique
Borel subalgebra with $N \in \mathfrak{b}$, and let $B$ be the Borel
subgroup corresponding to $\mathfrak{b}$.  Let $H$ be a Hessenberg space
for this Borel subalgebra.  The intersection of the Bruhat 
decomposition with respect to $B$ and the Hessenberg variety $\H(N,H)$ 
is a paving by affines of $\H(N,H)$ for each $H$.  
The nonempty cells of this paving are $B \pi B \cap \H(N,H)$ satisfying 
$\pi^{-1} \cdot N \in H$. 
\end{thmnonumber}

The dimension of each cell is the cardinality of a certain set of positive
roots depending on $\pi$, $H$, and $N$.  This set is described precisely in
Theorem \ref{reg nilps theorem}.
One consequence is that regular nilpotent, semisimple, and all type-$A_n$ 
Hessenberg varieties have no odd-dimensional cohomology.
In type $A_{n-1}$, I offer an alternative description of the paving
in terms of certain fillings of certain Young diagrams.  In this type,
the Hessenberg space $H$ is equivalent to a function $h: \{1,2,\ldots,n\}
\longrightarrow \{1,2,\ldots,n\}$ such that $h(i) \geq \max \{i,h(i-1)\}$
for all $i$.  (The relation between $h$ and $H$
 is described in greater detail in Section \ref{combi and Peterson}.)  
The theorem for nilpotent Hessenberg varieties follows.  

\begin{thmnonumber}
Let $N$ be a nilpotent operator.  Associate to $N$ the Young diagram
whose $i^{th}$ column has the same number of boxes as the dimension of
the $i^{th}$ Jordan block for $N$.  Assume this Young diagram is left-aligned
and bottom-aligned.

The Hessenberg variety $\H(N,H)$ is paved by affine spaces each of which
is associated to a permutation $\pi$.
The nonempty cells of the paving correspond to those fillings of the Young
diagram associated to $N$ for which the configuration
\[
\begin{tabular}{|c|}
\cline{1-1}  $\pi^{-1} k$  \\
\cline{1-1}  $\pi^{-1} j$  \\
\hline
\end{tabular}\]
only occurs if $\pi^{-1} j \leq h(\pi^{-1} k)$.

Given a nonempty cell represented as a (filled) Young tableau, the
dimension of this cell is the sum of the following two quantities:
\begin{enumerate}
\item The number of configurations
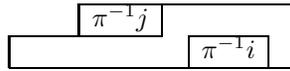
\[
\begin{tabular}{|cccccc|}
\cline{3-6}  \multicolumn{2}{c|}{  } & \multicolumn{1}{|c|}{$\pi^{-1} j$} & \multicolumn{2}{c}{   } & \\
\cline{1-3} \cline{5-5}  \multicolumn{2}{|c}{ $\hspace{1.6em}$ }& & &\multicolumn{1}{|c|}{$\pi^{-1} i$} &  \\
\hline
\end{tabular}\]
where box $i$ is to the right of or below box $j$, there is no box above $j$,
and the values filling these boxes satisfy $\pi^{-1} i > \pi^{-1} j$. 
\item The number of configurations
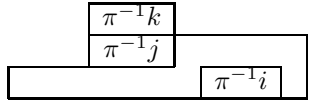
\[
\begin{tabular}{|cccccc|}
\cline{3-3}  \multicolumn{2}{c|}{$\hspace{2em}$  } & \multicolumn{1}{|c|}{$\pi^{-1} k$} & \multicolumn{3}{|c}{  } \\
\cline{3-6}  \multicolumn{2}{c|}{  } & \multicolumn{1}{|c|}{$\pi^{-1} j$} & \multicolumn{2}{c}{   } & \\
\cline{1-3} \cline{5-5}  \multicolumn{2}{|c}{ $\hspace{1.6em}$ }& & &\multicolumn{1}{|c|}{$\pi^{-1} i$} &  \\
\hline
\end{tabular}\]
where box $i$ is to the right of or below box $j$ and the values filling these
boxes satisfy $\pi^{-1} j < \pi^{-1} i \leq h(\pi^{-1} k)$.
\end{enumerate}
\end{thmnonumber}

This combinatorial 
method lends itself to computational results, as I demonstrate by providing
the Betti numbers of the Peterson variety in type $A_n$.  

The strategy of the proof is to use $M$ to choose a Bruhat 
decomposition so that on each Schubert cell the Hessenberg variety $\H(M,H)$
is an iterated tower of affine fibrations.  This procedure is independent
of the particular Hessenberg space $H$ so that if $M$ is fixed, 
the inclusion of Hessenberg spaces
gives rise to a natural inclusion of cells within their respective Hessenberg
varieties.  

Consider the nilpotent matrix
\[N = \left( \begin{array}{cccc} 0 & 1 & 0 & 0 \\
  0 & 0 & 1 & 0 \\ 0 & 0 & 0 & 1 \\ 0 & 0 & 0 & 0 \end{array} \right)\] and
the unipotent matrix 
\[u = \left( \begin{array}{cccc} 1 & a_{12} & a_{13} & 
a_{14} \\
0 & 1 & a_{23} & a_{24} \\ 0 & 0 & 1 & a_{34} \\ 0 & 0 & 0 & 1 
\end{array} \right).\]
Note that the conjugate $u^{-1} N u$ is
\[u^{-1} N u = \left( \begin{array}{cccc} 0 & 1 & a_{23}-a_{12} & 
a_{24}-a_{12}(a_{34}-a_{23})-a_{13}\\
0 & 0 & 1 & a_{34}-a_{23} 
\\ 0 & 0 & 0 & 0 \end{array} \right).\]  
Each flag $gB$ can be 
written as $u \pi B$ for some $u$ and a unique permutation $\pi$.  The flag
$gB$ is in the Hessenberg variety $\H(N,H)$ if and only if $u^{-1}Nu$ is
in $\pi H \pi^{-1}$.

Chapter \ref{defs and props} establishes the basic notational conventions
of this paper.  In type $A_{n-1}$,
we can identify the Lie algebra $\g{}$ with a subset of $n \times n$ matrices.
  This has a natural basis
of matrix units $E_{ij}$ defined to have value one in the $(i,j)$ entry
and zero elsewhere.  In terms of the previous example,
Chapter \ref{defs and props} shows that $H$ is spanned by
certain $E_{ij}$ and so the flag $gB$ is in $\H(N,H)$ if and only if
the matrix $u^{-1}Nu$ is zero in certain entries, which is equivalent to
certain polynomial equations in the entries of the matrix $u$ being zero.  
Section \ref{adjoint action of rows} describes
these equations in general.  These equations are
not necessarily linear, as the example shows.  However, the 
equations in the top row are affine functions of the variables $a_{1j}$ in 
terms of the variables $a_{ij}$ for $i \geq 2$.  Theorem \ref{seqlin} 
makes this claim in more general terms.  To prove Theorem \ref{seqlin}, we
need two main tools:
Section \ref{rows} describes a decomposition of
classical Lie algebras that generalizes the rows of a matrix; and
Section \ref{sequentially linear section} introduces a class of algebraic 
varieties called sequentially linear varieties whose added structure can be 
used to identify pavings.  
In Section \ref{criteria for affine}, we use the row decomposition
to show that Hessenberg varieties are paved by sequentially linear varieties 
and provide some
conditions under which they are in fact paved by affine spaces.  This 
amounts to partitioning flags $gB$ into Schubert cells and then showing that
within each Schubert cell, the affine function 
$a_{24}-a_{12}(a_{34}-a_{23})-a_{13}$ of the variables $a_{1j}$ will have
the same dimension solution space independent of the choice of $a_{ij}$ for
$i \geq 2$.  
Chapter 
\ref{main theorems section} contains the main theorems of this paper.  
Section \ref{combi and Peterson} includes a detailed analysis of 
the Peterson variety
in type $A_n$ and an explicit description of its cells.

\section{Definitions and Basic Properties}\label{defs and props}

Throughout this paper we use the notation and language of algebraic
groups as in \cite{H2}.

Let $G$ \label{G}be a linear algebraic group of classical type over the 
field $\mathbb{C}$ and
denote its Lie algebra by $\g{}$. \label{lie algebra}  
(The results in this paper for nilpotent
Hessenberg varieties generalize to fields of nonzero characteristic.
The results for semisimple Hessenberg varieties hold for
algebraically closed fields other than $\mathbb{C}$.)
%%Can use K char zero for much of this (applies to nilps, semisimps, not general)
Choose a maximal Cartan subalgebra 
$\mathfrak{h}$ 
\label{cartan subalgebra} 
in the Lie algebra and define positive roots
$\Phi^+$ \label{positive roots} 
and simple roots $\Delta$ \label{simple roots} 
with respect to this torus.
Write the decomposition of the Lie algebra into root spaces 
as $\g{} = \mathfrak{h} \oplus \bigoplus_{\alpha \in \Phi} \g{\alpha}$. 
\label{root space}
On occasion, we will fix a nonzero root vector $E_{\alpha}$ \label{root vector}
in $\g{}$
which spans $\g{\alpha}$.
Let $\mathfrak{b}$ \label{borel subalgebra} be the Borel subalgebra associated 
to $\Phi^+$ and let $\n{}$ \label{nilradical}
be its nilradical.  
%%
%%Note that the coordinate ring of $\n{}$ can be
%%written 
%%\[K[X_{\alpha}: \alpha \in \Phi^+].\]
%%This will be our convention throughout this paper.
%% *** check what these $X_{\alpha}$ actually are: characters?
%% co-characters?  eigenfunctions associated to $\g{\alpha}$?

Assume that the simple roots $\alpha_1, \ldots, 
\alpha_n$ are indexed according to the conventions of, e.g. \cite{H1}.
In other words, the bond between $\alpha_{n-1}$ and $\alpha_n$ in the
Dynkin diagram for $\g{}$ determine the type of the Lie algebra.

A Hessenberg space $H$ is a $\mathfrak{b}$-submodule of $\g{}$ which contains
$\mathfrak{b}$. \label{H} Let
\[\M_H = \{ \alpha \in \Phi: \g{\alpha} \subseteq H \}. \label{Hessenberg roots}\]
Then $\M_H$ is a subset of roots closed under addition of positive roots
and containing all positive roots.  Conversely, given any such subset of
roots, there is a unique Hessenberg space containing the corresponding root
spaces \cite{MPS}.  Consider the subspace $\mathfrak{i}$ orthogonal to $H$
 with respect to the Killing form in $\g{} = H \oplus \mathfrak{i}$.  
Note that $H$ is a Hessenberg space if and only if 
$\mathfrak{i}$ is an {\em ad}-nilpotent ideal in the sense
of \cite{CP}
with respect to the opposite Borel subalgebra $\mathfrak{b}^-$.  
The results of \cite{CP} thus show that the number of Hessenberg spaces
in type $A_{n-1}$ is the $n^{th}$
Catalan number.  They also show that slight
variations of Catalan numbers enumerate Hessenberg spaces 
in the other classical types.

Given a Hessenberg space $H$ and an element $M$ in $\g{}$, consider 
the subset of $G$ defined by
\[G(M,H) = \{g \in G: \Ad g^{-1}(M) \in H\}.\]
We often denote $\Ad g(M)$ by $g \cdot M = gMg^{-1}$.\label{g cdot M}

Since $H$ is closed under the adjoint action of the Borel subgroup $B$ 
corresponding to $\mathfrak{b}$, the subset $G(M,H)$ is closed under
right multiplication by $B$.  We may thus look at the image of $G(M,H)$
under the quotient map
\[\begin{array}{ccc} G &\supseteq &G(M,H) \\
                   \downarrows && \downarrows \\
                G/B & \supseteq& \H(M,H) \end{array} \]
The image $\H(M,H)$ is the Hessenberg variety corresponding to $M$
and $H$.  The space $G(M,H)$ is defined by closed conditions and the quotient
map is closed, so $\H(M,H)$ is a closed and hence projective variety.
%% Could cite closed quotient map: Massey, Alg. Top.: an Intro, p. 245

Unless otherwise stated, we assume that $M$ has been chosen from the Borel 
subalgebra $\mathfrak{b}$.  We use $N$ to denote an element of the
nilradical
$\bigoplus_{\alpha \in \Phi^+} \g{\alpha}$ and $S$ to denote an
element from the Cartan subalgebra $\mathfrak{h}$ in $\mathfrak{b}$.

\subsection{A Decomposition of the Nilradical
of $\mathfrak{b}$}\label{rows}

In this section we examine a decomposition of the nilradical 
$\mathfrak{n}$ of the fixed Borel algebra $\mathfrak{b}$ and prove some basic
properties of this decomposition.

The standard partial order on the set of roots is defined by
\[\alpha \geq \beta \mbox{  if and only if  } \alpha - \beta 
                \mbox{  is a sum of positive roots}. \label{root order}\]
We define $\alpha > \beta$ analogously so that $\alpha > \beta$
if and only if $\alpha \geq \beta$ and $\alpha \not = \beta$.
We often use the stronger condition that
$\alpha - \beta \in \Phi^+$.

Recall that the roots associated to the classical groups are 
described by the strings of simple roots given in this table.

\label{rootstructure}
\[ \begin{tabular}{|c|c|c|}
        \cline{1-3} Root & Parameters & Type \\
        \cline{1-3}       
        $ \sum_{j=i}^k \alpha_{j} $ &
        $1 \leq i \leq k \leq n$ &
        $A_n, B_n, C_n, D_n$ \\
        \cline{1-3}\multicolumn{3}{|c|}{\textup{except 
          $\alpha_{n-1}+\alpha_n \notin \Phi$ in type $D_n$}} \\
        \cline{1-3}
        $ \sum_{j=i}^n \alpha_{j} + \sum_{j=k}^n \alpha_j$ &
        $1 \leq i < k \leq n$ &
        $B_n$ \\
        \cline{1-3}
        $ \sum_{j=i}^n \alpha_{j} + \sum_{j=k}^{n-1} \alpha_j$ &
        $1 \leq i \leq k < n$ &
        $C_n$ \\
        \cline{1-3}
        $\sum_{j=i}^{n-2}\alpha_{j} + \alpha_n$ &
        $1 \leq i \leq n-2$ &
        $D_n$ \\        
        \cline{1-3}
        $ \sum_{j=i}^{n} \alpha_{j}  + 
        \sum_{j=k}^{n-2} \alpha_{j}$ &
        $1 \leq i < k \leq n-2$ &
        $D_n$ \\
        \hline
\end{tabular} \]

\noindent This follows from the definition of the 
root systems of classical Lie groups over characteristic zero fields as
in \cite[section 12]{H1}.

We often refer to the {\em extremal} simple roots of a root $\alpha$. 

\begin{definition}
The extremal simple roots of $\alpha$
are the simple roots $\alpha_i$ such that $\alpha - \alpha_i$ is
in $\Phi^+$. The extremal
roots of $\alpha$ are the positive roots $\beta$ such that $\alpha - \beta$
is in $\Phi^+$. 
\end{definition}

For instance, any non-simple root $\alpha$ of type $A_n$ has exactly
two extremal simple roots.  Recall that if $\alpha$ is
written as a sum of simple roots $\alpha = \sum_{j=1}^k \alpha_{i_j}$ 
then the height $\h (\alpha)$ of $\alpha$ is defined to be the number $k$ 
of simple summands.  
In type $A_n$, any non-simple root $\alpha$ 
has $2 \left(\h (\alpha) -1 \right)$ extremal roots.  
By inspection of Table \ref{rootstructure}, we see that in 
the other classical types a non-simple root can have either 
one, two, or three extremal simple roots.  

We define a partition of the positive roots and a collection of
nilpotent subalgebras associated to each part.
Let $\Phi^i$ be the subset of roots given by
\[ \Phi^i  = \{ \alpha \in \Phi^+: \alpha_i \leq \alpha,
                 \alpha_{j} \not < \alpha \textup{  for each  } j <i\}. 
\label{row roots}\]
In type $A_3$ this partition is $\Phi^1 = \{\alpha_1, 
\alpha_1+\alpha_2, \alpha_1 + \alpha_2 + \alpha_3\}$, 
$\Phi^2 = \{\alpha_2, \alpha_2 + \alpha_3\}$, and 
$\Phi^3 = \{\alpha_3\}$.  By contrast, the partition is
$\Phi^1 = \{\alpha_1, \alpha_1+\alpha_2, \alpha_1 + 2\alpha_2\}$
and $\Phi^2 = \{\alpha_2\}$ in type $B_2$. 

Let $\n{i} = \bigoplus_{\alpha \in \Phi^i} \g{\alpha}$ \label{lie algebra row}
be the subspace
of $\n{}$ spanned by the root spaces corresponding to $\Phi^i$.
Note that $\n{i}$ is a nilpotent subalgebra of $\n{}$.  
We refer to this subalgebra $\n{i}$ as the 
$i^{th}$ {\em row} of the Lie algebra
$\g{}$.  The terminology is inspired by the example of $\mathfrak{gl}_n$
considered as the collection of $n \times n$ matrices.  In this case, 
the subalgebra $\n{i}$ 
is precisely those matrices whose only
nonzero entries are in the $i^{th}$ row and above the diagonal.

The next lemma proves that
$\n{i}$ is either abelian or Heisenberg in classical types.

\begin{lemma}\label{rowstructure}
The subalgebras $\n{i}$ satisfy
\[ [\n{i},\n{j}] \subseteq \n{i} \textup{   for all } i \leq j.\]
\begin{enumerate}
\item In types $A_n$, $B_n$, and $D_n$, the $\n{i}$
        are abelian Lie algebras.
\item In type $C_n$, the $\n{i}$ are Heisenberg
        Lie algebras for $i < n$.  The subalgebra $\n{n}$
        is an abelian Lie algebra.
\end{enumerate}
\end{lemma}

\begin{proof}
The first claim follows from the definition of $\n{i}$ as well as the
property that
\begin{equation}\label{bracketrule}
 [\g{\alpha},\g{\beta}] = \left\{ \begin{array}{rl}
        \g{\alpha + \beta}  &\mbox{if}
        \hspace{1em} \alpha+\beta \in \Phi, \mbox{  and} \\
        0 &\mbox{if} \hspace{1em} \alpha + \beta \notin \Phi.
        \end{array} \right. \end{equation}
(See, e.g., \cite[section 8.4]{H1}.)

The second follows from the observation that 
\[\{\alpha+\beta: \alpha, \beta
\in \Phi^i; \alpha+\beta \in \Phi\} = \emptyset\]
in types $A_n$, $B_n$, and $D_n$.  In type $C_n$ the set 
\[\{\alpha+\beta: \alpha, \beta
\in \Phi^i; \alpha+\beta \in \Phi\} =
\left\{ \sum_{j=i}^{n-1} 2 \alpha_{j} + \alpha_n \right\} \subseteq \Phi^i.\]
Moreover, each root 
$ \alpha \not = \sum_{j=i}^{n-1} 2 \alpha_{j} 
+ \alpha_n$  in $\Phi^i$ generates a complementary root 
$\sum_{j=i}^{n-1} 2 \alpha_{j} + \alpha_n - \alpha$ in $\Phi^i$.  With
Property \eqref{bracketrule}, these conditions characterize abelian and
Heisenberg Lie algebras, respectively.
\end{proof}

In type $C_n$, the roots $2 \sum_{j=i}^{n-1} \alpha_j + \alpha_n$ are the
long roots.  We remark that there is exactly one long root in
each row $\Phi^i$ when $i<n$ in $C_n$.

The following proposition lists some characteristics of the row partition
in classical types.

\begin{proposition} \label{vertical rows}
Partition each row $\Phi^i$ of a given root system by height
and denote the parts
\[\Phi^i_k = \{\alpha \in \Phi^i: \h (\alpha) = k\}. \label{row height partition}\] 
These $\Phi^i_k$
satisfy the following properties in classical types:
\begin{enumerate}
\item If $\alpha$ is in $\Phi^i_k$ and $\beta$ is in $\Phi^i_{k-1}$
then $\alpha > \beta$.
\item For each $i$, the cardinality $| \Phi^i_k |$ is one except for at
most one $k_0$, for which $| \Phi^i_{k_0} | = 2$.
\item For all $2 \leq i \leq n-2$, if $| \Phi^i_k | = 2$ then 
$\Phi^{i-1}_{k+1} = \alpha_{i-1} + \Phi^{i}_{k}$ and
$|\Phi^{i-1}_{k+1}| = 2$.
\end{enumerate}
\end{proposition}

Root systems whose rows satisfy these properties are called {\em vertical}.  We
often call the rows themselves vertical.
To see that the rows in types $A_n$, $B_n$, $C_n$, and $D_n$ are 
all vertical, we simply inspect the entries in Table \ref{rootstructure}.  
Indeed, each row in types $A_n$, $B_n$, and $C_n$ is in fact 
ordered by height.  Conditions (2) and (3) apply only to rows in type $D_n$;
the conditions seem clumsy but will permit a general approach later in Lemma 
\ref{constantrank on cell}.

We often extend the definition of rows as follows to certain subgroups
in the unipotent subgroup 
of the 
linear algebraic group $G$ which corresponds to $\g{}$.  
When $G$ is of classical type other than $A_n$, 
we may assume that $G$ has been
embedded into $GL(N,\mathbb{C})$  
so that $\rk (G) = \lfloor N/2 \rfloor$
and so that the simple roots $\alpha_i$ for $G$ are simple
roots for $GL_N(\mathbb{C})$ when 
$i < \lfloor N/2 \rfloor$.
%% *** could reference that this is possible (G closed), Humphreys AG
Recall that 
$\exp (X) = \sum_{n \geq 0} \frac{X^n}{n!}$ is a formal power
series over $\mathfrak{gl}_N(\mathbb{C})$ which is a polynomial
whenever $X$ is nilpotent \cite[section 15.1]{H2}.
Define $U_i$ \label{group row} to be the subgroup generated by
$U_i = \exp \n{i}$. 
The map $\exp$ is a homomorphism when 
$\n{i}$ is an abelian Lie algebra.  Whether $\n{i}$ is abelian or
Heisenberg, the subgroup $U_i$ is the product of
the root subgroups $U_{\alpha}$ associated to the roots 
$\alpha$ in $\Phi^i$.
%% *** can state more by using 
%% the results in Humphreys, Alg. Gps. 15., which give you
%% root subgroups, and Chevalley commutator identity, which
%% tells you that each $U_i$ is in fact a group
Note that the rows $U_i$ generate the unipotent subgroup
$U = \prod_{i=1}^{n} U_{n-i+1}$.  We use this ordering
to describe $U$ throughout this paper.
%% *** check reference that the ordering by rows gives $U$

\subsection{The Bruhat Decomposition} \label{bruhat}

Here we recall some facts about Bruhat decompositions of the flag variety.
Write $T$ for the torus in $G$ whose Lie algebra is $\mathfrak{h}$ and 
denote the normalizer of $T$ by $N(T)$.  The Weyl group $W$ 
of $G$ is the quotient $W = N(T)/T$.
The Schubert cell in $G$ associated to a Weyl group
element $\pi$ is the double coset $B \pi B$.  By a slight abuse of
notation, we also denote the image of
this double coset under the projection to $G/B$ by $B \pi B$.
This is the Schubert cell corresponding to $\pi$ in the flag variety.

\begin{definition}
A paving ${\mathcal P}$ of an algebraic variety $X$ is an ordered partition
${\mathcal P} = ({\mathcal P}_1, {\mathcal P}_2, \ldots)$ of $X$ into 
disjoint varieties ${\mathcal P}_i$
such that
each finite union $\cup_{j\leq i} {\mathcal P}_j$ is closed in $X$.  If
each ${\mathcal P}_i$ is isomorphic to affine space, then ${\mathcal P}$
is a paving by affines.
\end{definition}

Pavings have less structure than CW-complexes but can still be used
to compute Betti numbers.  This motivates us to pave varieties by
simple spaces.

For instance, the Schubert cells $B \pi B$ form a paving by 
affines of the flag variety \cite[section 9.4]{Fu}.
Since Hessenberg varieties $\H(M,H)$ are closed in $G/B$ the
Schubert cells form a paving of $\H(M,H)$ as well.  The main claim 
of this paper is that in many cases $B$ can be chosen so that
this is in fact a paving by affines.

Define a subgroup $U_{\pi}$ \label{unipotent schubert cell} of the unipotent group $U$ by
\[ U_{\pi} = \{u \in U: \pi^{-1} \cdot u \in U^-\}\]
where $U^-$ is the opposite unipotent group associated to $U$.  The
group $U_{\pi}$
pa\-ra\-metrizes the Schubert cell corresponding to $\pi$
in the flag variety.  Note that
$U_{\pi}\pi$ is a set of coset representatives of the flags in 
the Schubert cell associated to $\pi$ \cite[sections 28.1 and 28.4]{H2}.
%%Can cite Cambridge text
Under the natural map
\[ \begin{array}{rl} U_{\pi} &\longrightarrow G/B \\
u & \mapsto u \pi B \end{array} \]
the subgroup $U_{\pi}$ is isomorphic to the corresponding Schubert
cell $B \pi B$ in the flag variety $G/B$ \cite[page 396]{FH}.
%% *** show holds for alg. grps. -- need char 0?
Denote the Lie algebra of $U_{\pi}$ by $\n{\pi}$.

\begin{proposition}
Each subgroup $U_{\pi}$ decomposes into a product of its rows
\[ U_{\pi} = \prod_{i=1}^n U_{\pi,n-i+1},\]
where $U_{\pi,i} = U_{\pi} \cap U_i$.
The Lie algebra $\n{\pi}$ can be written
\[\n{\pi} = \sp \langle \g{\alpha}: \alpha>0, 
        \pi^{-1} \alpha < 0  \rangle\]
and decomposes into rows $\n{\pi,i} = \n{\pi} \cap \n{i}$
each of which is the Lie algebra of the corresponding
subgroup $U_{\pi,i}$.
\end{proposition}

\begin{proof}
The subgroup $U$ can be written as a product
\[U = \prod_{i=1}^n \prod_{\alpha \in \Phi^{n-i+1}} \exp (\g{\alpha})\]
for this fixed ordering of $\Phi^+$ by rows.  This follows
from repeated application of the Chevalley commutator
relations.  Thus $U_{\pi}$
inherits a decomposition into row subgroups.

Moreover, 
\[\pi^{-1} \cdot \exp(\g{\alpha}) =  \exp(\g{\pi^{-1} \alpha}).\]
By the Chevalley commutator relations, the product 
$\prod \exp (X_{\alpha})$ is in $U_{\pi}$ if and only if 
each $X_{\alpha}$ is in $\n{\pi}$ \cite[section 26.3]{H2}.  
So $\n{\pi}$ is in fact the Lie algebra of $U_{\pi}$.
It follows that $\n{\pi,i}$ is the Lie algebra of
$U_{\pi,i}$ \cite[section 13.1]{H2}.
\end{proof}

\subsection{Sequentially Linear Varieties} \label{sequentially linear section}

In this section we define
 {\em sequentially linear} algebraic varieties and give some
of their preliminary properties.  We will later show that
Hessenberg varieties are examples of sequentially linear
varieties and use these properties to prove the main claims of 
this paper.

Let $X$ be an algebraic variety, either affine or projective.

\begin{definition}
A sequentially linear structure on a variety $X$ is a finite sequence
of varieties $X^i$ and morphisms $p_i$ so that 
\[X = X^n \stackrel{p_n}{\longrightarrow} X^{n-1}
   \stackrel{p_{n-1}}{\longrightarrow} \cdots X^1 
   \stackrel{p_1}{\longrightarrow} X^0 = \{ \textup{point}\}\]
and so that each $p_i$ has affine spaces as fibers.

If in addition each $p_i$ is a trivial affine fibration then
$X$ is a constant rank sequentially linear variety, often
simply called constant rank.
\end{definition}

The following proposition is clear from the definitions.

\begin{proposition} \label{seq linear with rank count}
If $X= X^n \stackrel{p_n}{\longrightarrow}
   \cdots X^1 
   \stackrel{p_1}{\longrightarrow} X^0 = \{ \textup{point}\}$ 
is a constant rank sequentially linear variety then $X$ is isomorphic
to affine $m$-dimensional space, where $m = \sum_{i=1}^n m_i$ and each
$m_i = \dim p_i^{-1} x_i$ for $x_i$ in $X^i$.
\end{proposition}
%%with associated
%%partitions $F$ and $\lambda$ and let $V$ be a closed set
%%on which $X$ has constant rank.  Suppose further that $I(V)$
%%is generated by $I(X)$ as well as a set of functions $F_V$ with 
%%\[F_V \subseteq \{f: f \in K[\lambda_i] \mbox{  for some  } 
%%                      i, \deg (f) = 1\}.\]  
%%Then $V$ is isomorphic to affine $m$-space where $m = n - \sum m_i$
%%and
%%\[m_i = \rk  \left\{ \begin{array}{c} \epsilon_x(F_i) = 0 \\
%%                              F_V \cap K[\lambda_i] = 0
%%      \end{array} \right\}\]
%%\end{proposition}
%% **** Comment: if $V = X$, then all conditions can be
%% satisfied except constant rank, but conclusions need not hold.
%% Also, if $I(X)$ is generated by $\{xy-1,z\}$, with two parts
%% and $\lambda_2 = \{x\}$, $F_2 = \{xy-1\}$, and $V$ is given
%% by $x-y$, then $F$ is constant rank on $V$, but $V$ is two
%% points and hence non-affine.  Thus, condition on $F_V$ needed. 
%%
%%\begin{proof}
%%Observe first that $m_i$ is an integer
%%since $\epsilon_x(F_i)$ has constant rank over $V$.
%%
%%Thus, we view $V$ as an iterated fibration with
%%$x_i$ chosen to satisfy
%%\[ \left\{ \begin{array}{c} \epsilon_x(F_i) = 0 \\
%%                              F_V \cap K[\lambda_i] = 0
%%      \end{array} \right\}\]
%%given solutions $x_1,\ldots , x_{i-1}$ to the previous
%%$i-1$ systems. By hypothesis, each system consists
%%of linear elements of $K[\lambda_i]$ so the solution space
%%to the $i^{th}$ system is affine.  By induction, $V$ is
%%an iterated fibration over an affine base space
%%with each fiber affine.  The proposition follows.
%%\end{proof}

\subsection{The Adjoint Action of Rows}\label{adjoint action of rows}

Here we discuss how the adjoint action $\Ad \! \! \!: G \longrightarrow \End{\g{}}$ 
behaves when considered as a map $\Ad \! \! \!: U_i \longrightarrow \End{\mathfrak{b}}$.
We also discuss the differential of this map $\ad \!: \n{i} \longrightarrow
\End{\mathfrak{b}}$.  A modification of the Chevalley commutator relation
and of Equation \eqref{bracketrule} permits an explicit description of $u^{-1}
\cdot M$ and $\ad X(M)$, respectively.
In both cases, properties of the $i^{th}$ row simplify this description
substantially.
Our ultimate goal is to use these
properties
to show that Hessenberg varieties are
paved by sequentially linear varieties.

As mentioned in Section \ref{rows}, the exponential map on a nilpotent
subalgebra $\n{}$ of $\g{}$ can be written $\exp(X) = 
\sum_{n \geq 0} \frac{X^n}{n!}$ (see, e.g.,
\cite[section 1.73]{K}).  
%% *** check category!
In particular, the operator $\ad X$ can be viewed
as an element of $\mathfrak{gl}(\mathfrak{b})$ and so $\exp \ad X$
is in $GL (\mathfrak{b})$.  We write this map explicitly as
\[\exp(\ad X) = \sum_{n \geq 0} \frac{(\ad X)^n}{n!}.\]

For any set $K$ of positive roots, the space $\n{K} = \oplus_{\alpha \in K}
\g{\alpha}$ is a vector subspace 
of $\mathfrak{b}$ whose natural basis of root vectors extends
to a basis for $\mathfrak{b}$.   Denote the corresponding 
quotient map by $\rho_{K}: \mathfrak{b} \rightarrows \n{K}$.
\label{projection to K} 
We may push $\rho_{K}$ forward to obtain the morphism 
\[\p{K}: \End(\mathfrak{b}) \longrightarrow \Hom(\mathfrak{b},\n{K}).\]
Here and subsequently 
$\End$ and $\Hom$ refer to the underlying 
vector-space endomorphisms and homomorphisms of the Lie algebras.  When 
$K = \Phi^i$ we abbreviate the projection to the $i^{th}$ row by $\rho_i$
\label{projection to row i}
and when $K = \{\alpha\}$ we write the projection to the root
space $\g{\alpha}$ by $\rho_{\alpha}$. \label{projection to alpha}
We also have occasion to write $\iota_K: \n{K} \hookrightarrow \mathfrak{b}$
\label{inclusion from K}
for the natural vector space inclusion.

\begin{lemma}\label{adform}
Let $X$ be an element of $\n{i}$.
In classical types the operator $(\ad X)^k $ in
$\End(\mathfrak{b})$ is identically zero when $k \geq 3$.  
When $k \geq 1$, 
\[\p{j} (\ad X)^k = 0 \textup{  for all  }j > i.\]
Furthermore,
\[ \p{i} (\ad X)^2 = 0 \textup{  in types $A_n$, $B_n$, and $D_n$}\]
and
\[ \textup{Im} \p{i} (\ad X)^2 \subseteq \g{\gamma_i} \textup{  in type $C_n$},\]
where $\gamma_i$ is the unique long root in $\Phi^i$. \label{long root row i}
\end{lemma}

\begin{proof}
Fix a set of generators $\{S_1, \ldots, S_{\scriptsize \rk G}\}$ 
for the torus $\mathfrak{h}$ in $\mathfrak{b}$.
The elements 
\[ \{(\ad X)^k E_{\alpha}: \alpha \in \Phi^+\}  \cup \{(\ad X)^k S_i : 1 \leq i \leq \rk G\}\]
generate the image 
$\mbox{Im} (\ad X)^k$.  
%%Given a subset $P \subseteq \Phi^+$, let $P^{i,k}$ be the 
%%subset of ordered $k$-tuples
%%\[
%%\begin{array}{ll}
%%P^{i,k} = \{(\beta_1, \ldots, \beta_k): &\beta_j \in \Phi^i  \hspace{.5em} 
%%      \forall j,
%%      \alpha + \sum_{j \leq l} \beta_j \in \Phi^+ \hspace{.5em} 
%%      \forall 1 \leq l \leq k \\ & \textup{  and  } \forall \alpha \in P\}.
%%\end{array}
%%\]
%%By abuse of notation, we permit $P = \emptyset$ and require in that
%%case the partial sums $\sum \beta_j$ to be in $\Phi^+$.
Using identity \eqref{bracketrule} repeatedly, we obtain
\begin{equation}\label{ad2}
(\ad X)^k E_{\alpha}  = \sum_{\scriptsize \begin{array}{c} \beta_1 + \cdots + 
\beta_k + \alpha \in \Phi \\
\beta_j \in \Phi^i \textup{  for all  } j \end{array}} 
        \left(\prod_{j=1}^k x_{\beta_j}\right) 
                        E_{\sum \beta_j + \alpha}.\end{equation}
Similarly,
\begin{equation}\label{adtorus}
(\ad X)^k S_i  = \sum_{\scriptsize \begin{array}{c}
\beta_1 + \cdots + \beta_k \in \Phi \\
\beta_j \in \Phi^i \textup{  for all  } j \end{array}}
        \left(\prod_{j=1}^k x_{\beta_j} \beta_j(S_i) \right) 
                        E_{\sum \beta_j}.\end{equation}
Table \ref{rootstructure} shows that no such
$\beta_1 + \cdots + \beta_k$ exist 
when $k$ is at least three, either for Equation \eqref{ad2} or for 
Equation \eqref{adtorus}.  It follows that
$(\ad X)^k$ is identically zero when $k \geq 3$.

Lemma \ref{rowstructure} and the definition of $\mathfrak{h}$ show that
$ [\n{i},\mathfrak{b}] \subseteq \sum_{j \leq i} \n{j}$.
Thus, when $X$ is in $\n{i}$ and $j>i$ the operator
 $\p{j} (\ad X)^k$ is identically zero.

In types $A_n$, $B_n$, and $D_n$ the Lie algebra $\n{i}$ is abelian so 
\[ \ad X \left( \sum_{j \leq i} \n{j} \right) \subseteq 
                \sum_{j < i} \n{j}\]
and $\p{i} (\ad X)^2$ is identically zero on $\mathfrak{b}$.

In type $C_n$ each row is a Heisenberg Lie algebra so 
\[ \ad X \left( \sum_{j \leq i} \n{j} \right) \subseteq 
                \sum_{j < i} \n{j} + \g{\gamma_i}. \]
\end{proof}

Several results follow.  The
following definitions are useful for notational brevity.

\begin{definition}
If $M$ is an element of $\g{}$ write 
$M = S_M + \sum_{\alpha \in \Phi_M} c_{\alpha} E_{\alpha}$
with the $c_{\alpha}$ nonzero constants and with $S_M$ in the fixed
Cartan subalgebra $\mathfrak{h}$.
The set $\Phi_M$ is the collection of roots associated to $M$.\label{roots of M}

Given a subset $\mathfrak{u} \subset \g{}$, write
$\Phi_{\mathfrak{u}} = \bigcup_{M \in \mathfrak{u}} \Phi_M$.
\label{roots of mathfrak u}
\end{definition}

In our applications $M$ is in $\mathfrak{b}$ and so
$\Phi_M$ is a subset of the positive roots.  We also write
$\Phi_{\pi}$ \label{roots of schubert cell}for 
$\Phi_{\n{\pi}}$, the roots associated to the
parameterization $U_{\pi}$ of the Schubert cell $B \pi B$.  
These roots are more concisely defined as $\Phi_{\pi} = \Phi^+ \cap 
\pi \Phi^- $ (see \cite[sections 28.1 and 28.4]{H2}).
Similarly $\Phi_{\pi,i} = \Phi_{\pi} \cap \Phi^i$ denotes
the roots associated to $\n{\pi,i}$.

\begin{corollary}\label{expadform}
Fix $X$ in $\n{i}$.  The operator
$\p{j} \exp \ad X$ in $\Hom(\mathfrak{b}, \n{j})$
satisfies the following:

\begin{enumerate}
\item If $j > i$ then
  $\p{j} \exp \ad X = \rho_j$. 
\item If $j=i$ then
\[ \p{i} \exp \ad X = \left\{ \begin{array}{ll}
                \rho_{i} + \p{i} \ad X 
                        &\hspace{1em}\mbox{ in types $A_n$, $B_n$, and $D_n$, and} \\
         \rho_{i} + \p{i} \ad X + \p{i} \frac{\scriptsize (\ad X)^2}{2} 
                        &\hspace{1em}\mbox{ in type $C_n$.}
        \end{array} \right.\]
\item If $j < i$ and 
$M = S_M + \sum_{\alpha \in \Phi_M} c_{\alpha}E_{\alpha}$
is in $\mathfrak{b}$ then
\[ \rho_{j} \exp \ad X (M) = \rho_{j} M + 
%%\ad X(S_M) + \frac{(\ad X)^2}{2}(S_M) \\
        \sum_{\scriptsize \begin{array}{c} \alpha \in \Phi_{M} \cap \Phi^j \\
               \ad E_{\alpha} (\n{i}) \not = \{0\} \end{array}}
                c_{\alpha} \left( 
                         \ad X(E_{\alpha}) + \frac{(\ad X)^2}{2}(E_{\alpha})
                         \right). \]
\end{enumerate}
\end{corollary}

\begin{proof}
Lemma \ref{adform} and the explicit description of
the exponentiation map show that 
\[\rho^*_j \exp \ad X = \left\{ \begin{array}{ll}
                \rho_j & \textup{if  }j > i, \\
                \rho_j + \rho^*_j \ad X & \textup{if  }j=i \textup{  in types } 
                     A_n,B_n,D_n, \textup{  and} \\
                \rho_j + \rho^*_j \ad X + \rho^*_j
                    \frac{\scriptsize (\ad X)^2}{2} & 
                    \textup{if  }j=i \textup{  in type  } C_n.
                \end{array} \right. \]
Equations \eqref{ad2} and \eqref{adtorus} complete the proof.
\end{proof}

\subsection{The Variety $U(M,\n{H})$}

In this section we show that Hessenberg varieties are paved by sequentially
linear varieties.  The strategy is
to intersect a fixed Hessenberg variety with
a fixed Schubert cell and study its preimage in $G$.  We then identify 
a subvariety
of the unipotent group in this preimage that is isomorphic to the original
intersection of Hessenberg variety with Schubert cell.  
This subvariety of $U$ will be sequentially linear.

%%Embed $\n{}$ in $K^{|\Phi^+|}$ using a fixed root space decomposition
%%so that 
%%\[K[\n{}] = K[X_{\alpha}: \alpha \in \Phi^+].\]
%%Write $\rho_{\alpha}: \n{} \rightarrows \g{\alpha}$ for the natural
%%projection to the root space corresponding to $\alpha$.
Fix $K \subseteq \Phi^+$ and define $\n{K}$ to 
be the subvariety of $\n{}$ given by $\n{K} = \bigoplus_{\alpha \in K}
\g{\alpha}$.  We define $\C{K} = \Phi^+ \setminus K$ to be the set
of positive roots complementary to $K$.\label{cal C: complement}
Let $M$ be an element of $\mathfrak{b}$ and write $U$ for the
unipotent subgroup corresponding to the nilradical $\n{}$.  Define
\[U(M,\n{K}) = \{u \in U: \Ad u^{-1} (M) \in \n{K} \}. \label{unipotent seq lin}\]
Recall that $\gamma_i$ denotes the longest root in $\Phi^i$ in type $C_n$.

\begin{theorem} \label{seqlin}
Fix $K \subseteq \Phi^+$.

The decomposition into rows defines a sequentially linear structure
on $U(M,\n{K})$ in types $A_n$, $B_n$, and $D_n$.  

In type $C_n$, suppose that whenever $\gamma_i \notin K$
and $\rho_{\gamma_i}(u^{-1} \cdot M) \not = 0$ 
for at least one $u$ in $U$ then either 
\begin{enumerate}
\item $\gamma_i(S_M) \not = 0$  or
\item both $\rho_{\alpha_i}(M) \not = 0$ and $(\gamma_i - \alpha_i)(S_M)=0$.
\end{enumerate}
Let $P^i_1 = \Phi^i - \{\gamma_i\}$ and $P^i_2 = \Phi^i - \{\gamma_i
- \alpha_i, \gamma_i\}$.  The refinement
of the decomposition into rows whose $2i+1^{th}$ part is $P^i_j$ and
whose $2i^{th}$ part is $\Phi^i-P^i_j$ defines a 
sequentially linear structure on
$U(M,\n{K})$ when condition (j) holds, for $j=1$ or
$j=2$.
\end{theorem}

\begin{proof}
Define 
$U^i = \prod_{j\geq i} U_j$ and $\n{K_i} = \left( \n{K} \cap \left(\bigoplus_{j
\geq i} \n{j}\right)\right) \oplus \bigoplus_{j<i} \n{j}$.  There is a natural
projection 
\[\begin{array}{ll}
p_i: & U^i\longrightarrow U^{i+1} \\
  &u_nu_{n-1} \cdots u_i \mapsto u_{n} \cdots u_{i+1} \end{array}\]
where elements of $U^i$ and $U^{i+1}$ are expressed as the ordered 
product of elements in decreasing rows.

Let $u$ be in $U^i$.  The calculations of Corollary \ref{expadform}
 show that
$\rho_j \left( u^{-1} \cdot M \right) = \rho_j \left( p_i(u)^{-1} \cdot M \right)$ 
for all $j>i$.  Thus, the
map $p_i$ restricts to a projection
\[p_i: U^i \cap U(M,\n{K_i}) \longrightarrow U^{i+1} \cap U(M,\n{K_{i+1}}).\]

We now inspect each fiber of this projection to ensure that it is an 
affine space.  Observe that $u$ is in
$U^i \cap U(M,\n{K_i})$ if and only if both $p_i(u)$ is in 
$U^{i+1} \cap U(M,\n{K_{i+1}})$ and $\rho_i(u^{-1} \cdot M) \in \n{K} \cap
\n{i}$.  Fix $u'$ in $U^{i+1}\cap U(M,\n{K_{i+1}})$
and write $M_i = (u')^{-1} \cdot M$.  Also write $u = u' u_i$ so
that $u^{-1}\cdot M = u_i^{-1} \cdot M_i$.  Note that there exists a 
unique $X_i \in \n{i}$ such that $u_i^{-1} = \exp X_i$ as described in Section
\ref{rows}.  Also note that 
\[u_i^{-1} \cdot M_i = \Ad \exp X_i (M_i) = \exp \ad X_i (M_i)\]
by \cite[section 1.93]{K}.

In types $A_n$, $B_n$,
and $D_n$, Corollary \ref{expadform} expands the projection 
$\rho_i(u_i^{-1} \cdot M_i)$ explicitly as the expression
\[\rho_i (u_i^{-1} \cdot M_i) = \rho_i \exp \ad X_i (M_i) = 
\rho_i M_i + \rho_i \ad X_i (M_i).\]  
Since $\ad X_i (M_i) = -\ad M_i (X_i)$, 
the second term is a linear function of $X_i$.  The first term simply
translates by the vector $\rho_i M_i$.  In other words, the set 
$\{u_i \in U_i: \rho_i (u_i^{-1} \cdot M_i) \in \n{K} \cap \n{i}\}$ 
describes an affine subspace of $U_i$ for each fixed 
 $u'$ in $U^{i+1} \cap U(M,\n{K_{i+1}})$.  This proves
the claim in those cases.

In type $C_n$,  define
\[V^i_j = \left( \prod_{\beta \in P^i_j}
   U_{\beta} \right) \cap U(M,\n{K_i} \oplus \n{\gamma_i}).\]
In each case we refine the tower of morphisms to include
\[U^i \cap U(M,\n{K_i}) \stackrel{p_i}{\longrightarrow} V^i_j
   \stackrel{p_{i,j}}{\longrightarrow} U^{i+1} \cap U(M,\n{K_{i+1}}).\]
Then $\rho_{\alpha}$ gives an affine transformation on $V_j^i$ of the
part of the $i^{th}$ row corresponding to $P^i_j$ 
for each $\alpha$ in $\Phi^+ - \{\gamma_i\}$.  Likewise, the function
$\rho_{\gamma_i}$ is an affine transformation in the entries corresponding
to $\Phi^{i} - P^i_j$ over the entries already fixed in $V^i_j$.  Both 
follow from Corollary \ref{expadform} and together prove the claim.
\end{proof}
%%
%%In each case, $F$ corresponds to a partition $\varphi = 
%%(\varphi_1, \ldots, \varphi_k)$ of $P_H$, where $\varphi_i$ is
%%defined by
%%\[F_i = \{ \rho_{\alpha} (u^{-1} \cdot M): \alpha \in \varphi_i\}.\]

The following lemma relates the varieties $U(M,\n{H})$ to Hessenberg
varieties.

\begin{lemma} \label{uni is hess}
Let $ \H(M,H) \cap B \pi B $ be the intersection of the Schubert
cell corresponding to $\pi$ with the Hessenberg variety $\H(M,H)$.  
Then
\[ \H(M,H) \cap B \pi B \cong U_{\pi} \cap U(M,\n{\pi \cdot H})\]
for
$\n{\pi \cdot H} =  \n{} \cap (\pi \cdot H)$. 
Consequently, the intersection of each Hessenberg variety with
each Schubert cell is sequentially linear in the classical types.
\end{lemma}

\begin{proof}
The Schubert cell $B \pi B$ is isomorphic to $U_{\pi}$ as 
discussed in \cite[page 396]{FH}.
Since the projection from $G(M,H)$ to the Hessenberg variety
$\H(M,H)$ is $\Ad(B)$-invariant, this isomorphism restricts to
the intersection
\[ \H(M,H) \cap B \pi B \cong U_{\pi} \cap U(M,\n{\pi \cdot H}).\]
Explicitly, if we write $gB$ for the flag corresponding to $g$ we
see that
\[\begin{array}{rlr}
        gB & \in \H(M,H) &\Leftrightarrow \\
        (gb)^{-1} \cdot M & \in H \hspace{1em} \textup{  for all  } b \in B & \Leftrightarrow \\
        (u \pi)^{-1} \cdot M & \in H \mbox{  for  } u \pi B = g B, u \in U_{\pi}
                                & \Leftrightarrow \\
        u &\in U(M,\n{\pi \cdot H}) \cap U_{\pi}. & \end{array}\]

The conclusion follows from Theorem \ref{seqlin} and the
definition of the Hessenberg space $H$.
\end{proof}

The following theorem summarizes these results.

\begin{theorem} \label{affines}
Suppose there exists a Borel subgroup $B$ such that
\[  U_{\pi} \cap U(M,\n{\pi \cdot H}) \cong
                \left\{ \begin{array}{l} \emptyset \\
                        \mathbb{C}^d \mbox{  for some $d$}
                        \end{array} \right.\]
for each $\pi$ in $W$.  Then the paving of
$\H(M,H)$ obtained by intersecting the Hessenberg variety 
with Schubert cells
is a paving by affines
such that the cell $\H(M,H) \cap B \pi B$ has dimension $d$.
\end{theorem}

\begin{proof}
The Schubert cells $\H(M,H) \cap B \pi B$ pave each
Hessenberg variety as per the comments in Section \ref{bruhat}.
Also, $\H(M,H) \cap B \pi B$ is isomorphic to
$U_{\pi} \cap U(M,\n{\pi \cdot H})$ by Lemma 
\ref{uni is hess}.  If the hypotheses hold
then the Bruhat decomposition actually
gives a paving by affines of $\H(M,H)$ and the dimension
of $\H(M,H) \cap B \pi B$ equals that of 
$U_{\pi} \cap U(M,\n{\pi \cdot H})$ for each $\pi$.
\end{proof}

For convenience, we remark that 
the Adjoint action of $G$ gives an action of $G$ on Hessenberg varieties
defined by $\Ad g^{-1} \left( \H(M,H)\right) = \H(g^{-1} \cdot M, g^{-1}
\cdot H)$.  
The action $\Ad g^{-1}$ is an isomorphism of Hessenberg varieties and so
all Hessenberg varieties in a fixed $G$-orbit are isomorphic to each
other.  We state
this as a lemma though the proof is immediate.

\begin{lemma}\label{Hess invariant of orbit}
If $g \in G$ then $\H(M,H) \cong \H(g^{-1} \cdot M, g^{-1}
\cdot H)$.
\end{lemma}

Note that $\H(M,H)$ is not independent of the choice of Borel 
subalgebra $\mathfrak{b} \subseteq H$.  Indeed, the definition of 
the Hessenberg space $H$ requires that $\ad 
\mathfrak{b} (H) \subseteq H$.  This is not generally true of Borel 
subalgebras contained in $H$.

We interpret the choice of a Borel $B$ in Theorem 
\ref{affines} as fixing a basis for the flag variety $G/B$ with
respect to which we
consider $\H(M,H)$.  We will use this in Chapter \ref{main theorems section} 
to select
a computationally convenient form of $M$ from its $G$-orbit.  

\subsection{Criteria for $U(M,\n{\pi \cdot H}) \cap U_{\pi}$ to be an affine space}\label{criteria for affine}

With certain extra assumptions on $M$, the variety $U(M,\n{\pi \cdot H})$
will be not just sequentially linear but will also intersect the
closed subgroup $U_{\pi}$ in an affine space.  
Theorem \ref{affines} will
then imply that $\H(M,H)$ is paved by affines whose dimensions we can 
identify.

Let $M$ be an element of $\mathfrak{b}$ written
\[M = S_M + \sum_{\beta \in \Phi_M} m_{\beta} E_{\beta} = S_M + N\]
for nonzero constants $m_{\beta}$, a semisimple element $S_M$ in
$\mathfrak{h}$, and a nilpotent $N$ in the nilradical $\n{}$.  

\begin{definition}
A collection of roots $P$ is non-overlapping if for
no pair $\alpha, \beta$ in $P$ is $\alpha > \beta$.

If $M=S_M + N$ is written as above 
then $M$ is non-overlapping if both of the following hold:
\begin{enumerate}
\item $\Phi_N$ is non-overlapping.
\item For each $\alpha \in \Phi_N$ and each simple root $\alpha_i$ with 
$\alpha \geq \alpha_i$ the equality $\alpha_i(S_M)=0$ holds.
\end{enumerate}
\end{definition}

By an abuse of notation, we call the roots $\Phi_M$ non-overlapping
if $M$ is non-overlapping.  
Note that the second condition implies that $\beta(S_M) = 0$
for each $\beta \leq \alpha$ and each $\alpha \in \Phi_N$.  Consequently,
the Lie algebra elements $S_M$ and $N$ commute.  However, the requirement that 
$\ad S_M (N)=0$ is not sufficient to ensure that the second condition holds.
For instance, in $\mathfrak{gl}_3$ the element 
\[M = \left( \begin{array}{ccc} x & 0 & 1 \\
0 & y & 0 \\ 0 & 0 & x \end{array} \right)\] is {\em not}
non-overlapping while \[M' = \left( \begin{array}{ccc} x & 1 & 0 \\
0 & x & 0 \\ 0 & 0 & y \end{array} \right)\] is non-overlapping.  Not
only are $M$ and $M'$ in the same $G$-orbit but they both have a Jordan
decomposition into diagonal and non-diagonal parts.  Thus, if 
$M$ is non-overlapping then $M=S_M +N$ as above is a Jordan decomposition
but not vice-versa.

\begin{lemma}\label{nonoverlap gives rootsofN}
Let $M$ be in $\mathfrak{b}$ and $u$ be in $U$.  If $\Phi_M$ is 
a non-overlapping set of roots then 
$\Phi_M \subseteq \Phi_{u^{-1} \cdot M}$ and
$\Phi_M \cap \Phi_{u^{-1} \cdot M - M}$ is empty.
\end{lemma}

\begin{proof}
Let $c_{\alpha}$ be nonzero constants so that 
\[u^{-1} \cdot M = S_M + \sum_{\alpha \in \Phi_{u^{-1} \cdot M}} 
        c_{\alpha} E_{\alpha}.\] 
%%Let $P_{\alpha}$ be the collection of ordered sets of roots defined by
%%\[P_{\alpha} = \{ (\beta_0, \ldots, \beta_j): j>0, \beta_0 \in \Phi_M,
%%              \beta_i \in \Phi_u \hspace{0.5em}\forall i>0,  \sum_i \beta_i =
%%               \alpha \}\]
%%and let $P_{\alpha}'$ be the collection of ordered sets of roots defined by
%%\[P_{\alpha}' = \{ (\beta_1, \ldots, \beta_j): j>0,
%%              \beta_i \in \Phi_u \hspace{0.5em}\forall i>0,  \sum_i \beta_i =
%%               \alpha \}.\]
%%Then for each root $\alpha$ in $\Phi_{u^{-1} \cdot M}$ there
%%exist subsets $C_{\alpha} \subseteq P_{\alpha}$ and 
%%$C_{\alpha}' \subseteq P_{\alpha}'$ such that
%%the coefficient of $E_{\alpha}$ in $u^{-1} \cdot M$ is
%%\[c_{\alpha} = m_{\alpha} + \sum_{ (\beta_i) \in C_{\alpha}' }
%%                      \prod_{i>0} u_{\beta_i} \beta_i(S_M)
%%                   + \sum_{ (\beta_i) \in C_{\alpha} } m_{\beta_0}
%%                      \prod_{i>0} u_{\beta_i}.\]
%%The set $P_{\alpha}$ is empty when $\alpha$ is in $\Phi_M$.  The 
%%product $\prod_{i>0}  \beta_i(S_M)$
%%is zero when  $(\beta_i)$ is in $C_{\alpha}'$.  
%%(Both facts follow from the definition of non-overlapping.)
Fix $\alpha$ in $\Phi_M$.  Write $u = u_n u_{n-1} \cdots u_1$ for each
$u_i \in U_i$ and apply the conclusions of Corollary \ref{expadform} 
repeatedly to 
$u^{-1} \cdot M = (\prod_{i=1}^n \exp X_{n-i+1})^{-1} \cdot M$.
For our purposes, we need only the result that 
the coefficient $c_{\alpha}$ is the sum of the following 
three quantities.  The first
is $m_{\alpha}$.  The second is a sum of terms of the form
$\prod_{i=1}^r u_{\beta_i}\beta_i(S_M)$ for $r$-tuples of
$\beta_i \in \Phi^+$ such that $\beta_1 + \cdots + \beta_r = \alpha$
and constant coefficients $u_{\beta_i}$.
Regardless of $r$ or the choice of $\beta_i$, 
Condition (2) in the definition of non-overlapping
ensures that this quantity is zero.
The third quantity is a sum of terms of the form
$m_{\beta_0} \prod_{i=1}^r u_{\beta_i}$ for $r$-tuples of $\beta_i \in
\Phi^+$ such that $\beta_0 + \beta_1 + \cdots + \beta_r = \alpha$ and
constants $u_{\beta_i}$.
Again by definition of non-overlapping, this quantity is zero.
Consequently the term $c_{\alpha}$ is simply $m_{\alpha}$.
\end{proof}

\begin{corollary} \label{empty cells}
If $\Phi_M$ is non-overlapping and $\Phi_M \not \subseteq K$ then 
\[U(M,\n{K}) = \emptyset.\]
\end{corollary}

\begin{proof}
If $u$ is in $U$ then $\Phi_M \subseteq \Phi_{u^{-1} \cdot M}$ 
by the previous lemma.  Since $\Phi_M \not \subseteq K$ the
element $u^{-1} \cdot M$ cannot be in $\n{K}$.
\end{proof}

With certain hypotheses we can reduce to the study of nilpotent
Hessenberg varieties.  For each semisimple element $S \in \mathfrak{h}$, 
write 
$\Phi_S^+$ for the set of roots
$\Phi_S^+ = \{\alpha \in \Phi^+: \alpha(S) = 0\}$. \label{roots annihilating S}  
Let $\Delta_j$ be the maximal
irreducible subsets of $\Delta \cap \Phi_S^+$.  Denote the parabolic
subalgebra
associated to the simple roots $\Delta_j$ by $\mathfrak{p}_{\Delta_j}$ and
choose its associated Levi part $\mathfrak{l}_{\Delta_j}$ so that
$\mathfrak{l}_{\Delta_j} \supseteq \mathfrak{h}$. 
Let $U_{\Delta_j}$ be
$L_{\Delta_j} \cap U$.
Write $\mathfrak{p}_M$ for the parabolic
subalgebra associated to $\bigcup \Delta_j$ and write $\mathfrak{l}_M$ and
$\n{M}$ for its associated Levi and nilpotent parts.  Recall that 
$\C{\pi \cdot H}$ is the set of positive roots complementary to $\pi \M_H
\cap \Phi^+$.

\begin{theorem} \label{reduce to nilp}
Let $M = \sum M_j$ be an element of $\mathfrak{b}$ and write the 
decomposition of each summand as $M_j = S_j + N_j$ for $S_j$ in 
$\mathfrak{h} \cap \mathfrak{l}_{\Delta_j}$ and $N_j$ in $\n{}$.  
Let $S = \sum S_j$
be the semisimple part of $M$.  Assume that
the following conditions hold:
\begin{enumerate}
\item The set of roots $\Phi_S^+$ is the union 
$\bigcup \sp \langle \Delta_j \rangle$, each $N_j$ is in the corresponding
$\mathfrak{l}_{\Delta_j}$, and $\Phi_{M_j}$ is non-overlapping.
\item The variety $U(N_j, \mathfrak{l}_{\Delta_j} \cap \n{\pi \cdot H}) \cap 
U_{\Delta_j} \cap U_{\pi}$ 
has a sequentially linear structure which refines 
the decomposition into rows of Theorem \ref{seqlin} and which is
constant rank of total dimension $m_j$.
\end{enumerate}
Then the variety $U(M, \n{\pi \cdot H}) \cap U_{\pi}$ is nonempty if and
only if $\pi^{-1} \cdot N_j$ is in $H$ for each $j$.  If nonempty, it is
an affine space of
dimension $|\Phi_{\n{M}} \cap \Phi_{\pi} \cap \pi \M_H| + \sum m_j$.
\end{theorem}

\begin{proof}
Choose any $u$ in $U$.
Consider the operator 
\[\rho^*_{\n{M} \cap \n{i}} 
\iota_{* \n{M} \cap \n{i}} \ad (u^{-1} \cdot M)
\in \End(\n{M} \cap \n{i})\]
written with respect to the basis of root
vectors in $\n{i}$ ordered by height.   In type $D_n$, fix any order 
among root vectors in $\n{i}$
of the same height.  
In type $C_n$, use the 
refinement of the row decomposition described in Case (1) of Theorem
\ref{seqlin}. 
We saw that this operator is an affine transformation 
in Theorem \ref{seqlin}.  Here we show that 
its solution space has the same dimension for any $u \in U_{\pi}$.

The operator $\ad (u^{-1} \cdot M)$ acts as dilation by the nonzero constant 
$\alpha(S)$ on each root vector $E_{\alpha}$ in the nilpotent 
subalgebra $\n{M}$.  This follows since the $\Delta_j$ generate both the 
roots associated to $\mathfrak{l}_M$ and $\Phi_S^+$. Moreover, the operator 
$\rho^*_{\alpha} \ad (u^{-1} \cdot M)$ is identically zero on $\g{\beta}$ for
each root $\beta \not \! \leq \alpha$ in $\Phi^i$.  
This shows that the operator $\rho^*_{\n{M} \cap \n{i}} 
\iota_{* \n{M} \cap \n{i}} \ad (u^{-1} \cdot M)$
is a lower triangular matrix with
nonzero entries $\alpha(S)$ along the diagonal with respect to the
basis defined above.  

It follows that
$U(M,\n{\pi \cdot H}) \cap \exp \n{M}$ is a constant rank sequentially
linear variety, which is to say an affine space of dimension
\[|\Phi_{\n{M}} \cap \Phi_{\pi}| - 
      |\Phi_{\n{M}} \cap \Phi_{\pi} \cap \C{\pi \cdot H}| = 
      |\Phi_{\n{M}} \cap \Phi_{\pi} \cap \pi \M_H|. \]

If $E_{\alpha}$ is in $\mathfrak{l}_{\Delta_j} \cap \n{}$ then $\alpha$
is generated by $\Delta_j$
by Equations \eqref{ad2} and \eqref{adtorus}.
Thus, the hypotheses on $S_j$ imply that
the operators $\rho^*_{\mathfrak{l}_{\Delta_j}}  
\ad (u^{-1} \cdot M_j)$
and $\rho^*_{\mathfrak{l}_{\Delta_j}} 
\ad (u^{-1} \cdot N_j)$ agree on the root space $\g{\alpha}$.
Consequently, the variety $U(N_j, \mathfrak{l}_{\Delta_j} 
\cap \n{\pi \cdot H}) \cap U_{\Delta_j} \cap U_{\pi}$ is a constant rank
sequentially linear variety  
if and only if the variety $U(M_j, \mathfrak{l}_{\Delta_j} 
\cap \n{\pi \cdot H}) \cap U_{\Delta_j} \cap U_{\pi}$ is.  
If so, their dimensions
are the same.

Furthermore, 
\[\rho^*_{\mathfrak{l}_{\Delta_j}} 
    \ad \left( u^{-1} \cdot \left(\sum_j N_j \right) \right) = 
    \rho^*_{\mathfrak{l}_{\Delta_j}} 
    \ad \left( 
\rho_{\mathfrak{l}_{\Delta_j}} u \right)^{-1} \cdot N_j\]
as operators in $\Hom(\n{}, \mathfrak{l}_{\Delta_j} \cap \n{})$.
It follows that $U(M,\mathfrak{l}_{\Delta_j} \cap \n{\pi \cdot H}) \cap 
U_{\Delta_j} \cap U_{\pi}$ 
is a constant rank sequentially linear variety 
if and only if each $U(N_j,\mathfrak{l}_{\Delta_j} \cap 
\n{\pi \cdot H}) \cap U_{\Delta_j} \cap U_{\pi}$ is.  
The dimension of the former is
the sum of the dimensions of the latter.

Take the sequentially linear structure on $U(M,\n{\pi \cdot H}) \cap
U_{\pi}$ obtained by refining the row decomposition of Theorem \ref{seqlin}
to first $\mathfrak{l}_M \cap \n{i}$ and then $\n{M} \cap \n{i}$.  The
arguments above show that this variety is constant rank.  It follows that
the variety is an affine space.
The total dimension of $U(M,\n{\pi \cdot H}) \cap U_{\pi}$ 
is obtained by summing the dimensions restricted to $\n{M}$ and each 
$\mathfrak{l}_{\Delta_j} \cap \n{}$.
\end{proof}

We now establish conditions under which a nilpotent Hessenberg
variety intersects Schubert cells in affine spaces.

%%\begin{definition}
%%A set of roots $P$ is saturated
%%in the set $P' \subseteq \Phi^+$ if there exists a
%%well-defined function $f: P' \rightarrow P$ such that
%%$f(\alpha)$ is an extremal root for $\alpha$.
%%\end{definition}
%%
%%Lemma \ref{nonoverlap gives rootsofN} implies that if 
%%$\Phi_N$ is non-overlapping it cannot be
%%saturated in the set $\Phi_{u^{-1} \cdot N}$.  By contrast, the set of 
%%simple roots $\Delta$ is saturated in any subset of $\Phi^+ - \Delta$.
%%
\begin{lemma} \label{sat gives rankcount}
Suppose that $P$ and $P'$ are subsets of $\Phi^+$
and suppose that each root $\alpha$ in $P'$ has an extremal root 
$\beta$ in $P$.  

Let $\pi$ be a Weyl group element and let $H$ be a Hessenberg space
in $\g{}$ with roots $\M_H$ such that $\pi P \subseteq \M_H$.  

If $\alpha$ is a root in $P'$ such that
\[\pi \alpha \in \Phi - \M_H\]
then there exists a root $\beta$ in $P$ satisfying
\[\alpha - \beta \in \Phi_{\pi}.\]
\end{lemma}

\begin{proof}
If $\beta$ is an extremal root for $\alpha$ then $\alpha - \beta$ is 
a positive root.  Let $\beta$ in $P$ be extremal for $\alpha$.  Note
that
$ \pi \alpha = \pi(\alpha - \beta) + \pi \beta $
is {\em not} in $\M_H$ by hypothesis.  Since
$\pi \beta$
is in $\M_H$ and since $\M_H$ is closed under addition by positive roots, 
the root $\pi (\alpha-\beta)$
must be negative.  This means that $\alpha-\beta$ is in $\Phi_{\pi}$.  
\end{proof}

We are building to a lemma that describes one set of conditions under which
the variety $U(N,\n{\pi \cdot H}) \cap U_{\pi}$ is an affine space 
when $N$ is nilpotent.  
%%
%%No longer have \label{constantrank on rowofcell}

Given a set $P \subseteq \Phi^+$ define
\[P^i_{-N} = \{\alpha - \beta: \alpha \in P, \beta \in \Phi_N, 
    \alpha - \beta \in \Phi^i\} \label{row i pm roots of N}\]
and define $P^i_{+N}$ to be the roots $\alpha+\beta$ satisfying the
analogous conditions.  We often suppress the superscript and 
write $P_{\pm N}$ if $P$ is a subset of $\Phi^i$.

\begin{lemma} \label{constantrank on cell}
Let $\Phi_N$ be a non-overlapping set of roots.
Assume that $\Phi^i_{k,-N}$ is either empty or $\Phi^i_j$ for some $j$.
When $| \Phi^i_k | = 2$, also assume that $\Phi^i_k = \Phi^i_{k+1,-N}
= \Phi^{i-1, i}_{k,+N}$.

Suppose that $\{\alpha\}_{-N}$ is nonempty for each $\alpha$ 
in $\Phi^i \cap \Phi_{U_{\pi} \cdot N} \cap \C{\pi \cdot H}$
 and that if $\alpha = \sum_{\beta_j \in \Phi^i} \beta_j$ then
each $\beta_j \leq \{\alpha\}_{-N}$.
Then $U(N, \n{\pi \cdot H}) \cap U_{\pi}$ is an affine space of
dimension
\[|\Phi_{\pi}| - |\C{\pi \cdot H} \cap \Phi_{U_{\pi} \cdot N}|.\]
\end{lemma}

\begin{proof}
Let $P^i_{\pi \cdot H,-N}$ denote the set $\left( \Phi^i \cap \C{\pi \cdot H}
\right)^i_{-N}$.
Consider the map
\[N' = \left( \rho^*_{\Phi^i \cap \Phi_{U_{\pi \cdot N}} 
\cap \C{\pi \cdot H}} \right) \left( 
\iota_{* P^i_{\pi \cdot H,-N}} \right) \ad (u^{-1} \cdot N) \in
\Hom(\n{P^i_{\pi \cdot H,-N}}, \n{\Phi^i \cap \Phi_{U_{\pi \cdot N}} 
\cap \C{\pi \cdot H}})\]
as it acts with respect to a basis
of root vectors ordered by height.  We will write the matrix for $N'$
explicitly and show that its rank is independent of the choice of $u$ so long
as $u \in U_{\pi}$.  In Theorem \ref{seqlin} we gave a sequentially linear
structure for $U(N,\n{\pi \cdot H})$ by showing that certain projections
$\rho_i(u_i^{-1} \cdot N_i)$ were affine transformations for each $i$.  Since
$N'$ is the linear part of this affine transformation, its rank is
constant on $U(N,\n{\pi \cdot H}) \cap U_{\pi}$ if and only if
$U(N,\n{\pi \cdot H}) \cap U_{\pi}$ is a constant rank sequentially linear
variety, i.e. an affine space.

First we characterize the matrix for $N'$.  
Order the rows by height according to $\Phi^i_{k}$,
fixing an order if $|\Phi^i_{k}| = 2$.  Note that if $\alpha$ is in 
$\Phi^i \cap \C{\pi \cdot H}$ then there exists $\alpha'$ in $\Phi_N$
such that $\alpha - \alpha'$ is in
$\Phi_{\pi}$ by Lemma \ref{sat gives rankcount}.
Consequently, the columns in the matrix are indexed by elements of
$\Phi^i_{-N}$.
Order the columns by increasing
$k$ in $\Phi^i_{k,-N}$ using the order inherited from the associated 
rows of the matrix if $|\Phi^i_k| =2$.  

We begin by showing that the order on the columns
is consistent with ordering the columns by height.
Suppose $\alpha$ is in $\Phi^i_k$ and $\beta$ is in $\Phi^i_j$, that
$\alpha > \beta$, and that neither $\{\alpha\}_{-N}$ nor
$\{\beta\}_{-N}$ is empty.  By the
verticality of the rows, we know each element of $\Phi^i_k$ is 
greater than each element of $\Phi^i_j$.  Since 
neither $\Phi^i_{k,-N}$ nor $\Phi^i_{j,-N}$ is empty they 
equal $\Phi^i_{k'}$ and $\Phi^i_{j'}$ respectively.  Let $\alpha'$ and
$\beta'$ be in $\Phi_N$ with the property that 
$\alpha - \alpha' \in \Phi^i_{k'}$ and $\beta - \beta' \in \Phi^i_{j'}$.  
The sets $\Phi^i_{k'}$ and $\Phi^i_{j'}$ are comparable by 
definition of verticality.  Moreover,
\[(\alpha - \alpha') - (\beta - \beta')  = (\alpha - \beta) + (\beta' - 
\alpha') >0\]
since $\Phi_N$ is non-overlapping.  Thus each element of $\Phi^i_{k,-N}$
is greater than $\Phi^i_{j,-N}$, which shows that the order on the columns
previously defined is in fact the order by height (when $|\Phi^i_k| \not = 2$).

There is no $\beta < \alpha'$ with $\beta$ in $\Phi_{u \cdot N}$ since 
that would imply there were a $\beta'$ in $\Phi_N$ with $\beta'<\beta<\alpha'$,
which contradicts the definition of non-overlapping.  
The operator $\rho_{* \alpha}  \ad (u^{-1} \cdot N)$ scales $\g{\alpha - 
\alpha'}$ by $n_{\alpha'} \not = 0$ as shown in Lemma 
\ref{nonoverlap gives rootsofN}.  
In sum, 
\[ \rho^*_{i}
\iota_{* \alpha - \beta} \ad (u^{-1} \cdot N) = 
\left\{ \begin{array}{ll} 0 & \beta < \alpha' \textup{  and }\\
                         n_{\alpha'} & \beta = \alpha'. \end{array} \right.\]
It follows that the matrix for $N'$ is lower triangular with respect to 
this basis.  

If no $|\Phi^i_k| =2$ then this matrix is nonzero along its 
diagonal.  This proves the claim in root systems for which 
 $\exp \ad (u^{-1} \cdot N) X_i$ 
is an affine transformation and no $|\Phi^i_k|=2$, namely in 
types $A_n$ and $B_n$.

In type $C_n$, the operator $\rho^*_{i} \exp \ad (u^{-1} \cdot N)$ is 
not an affine transformation.  However, the projection $\rho^*_{\alpha}
\exp \ad (u^{-1} \cdot N)$ is affine for each $\alpha \not = \gamma_i$.
Moreover, the non-affine function $(\ad X)^2$ satisfies $\rho^*_i (\ad X)^2 =
\rho_i^* (\ad (\rho_{\scriptsize \beta: \h \beta < \h \alpha'} X))$ for 
all $X$ in $\n{i}$.  In other words, the map $\rho^*_{\gamma_i} \iota_{* \{
\gamma_i-\alpha_i, \gamma_i\}} \exp \ad (u^{-1} \cdot M)$ is linear and 
a matrix of
constant rank.  Using the sequentially linear 
structure of Theorem \ref{seqlin},
the claim follows in type $C_n$, as well.

When there exists $\Phi^i_k=\{\beta_1, \beta_2\}$,  
note that 
\[\left( \frac{1}{n_{\beta_1}}\rho_{* \beta_1} - 
\frac{1}{n_{\beta_2}} \rho_{* \beta_2} \right) \ad (u^{-1} \cdot N) \]
is an affine transformation of the root vectors corresponding to 
$\Phi^{i+1}_{k-1}$ that depends on the choices for the root vectors 
corresponding to $\bigcup_{j<k-1} \Phi^i_j \cup \bigcup_{j>i} \Phi^j$.
This transformation must be identically zero if $\Phi^i_k \subseteq \C{\pi 
\cdot H}$.
Moreover, this affine transformation
is linearly independent from $\rho_{* \alpha} \ad u^{-1} \cdot N$ for $\alpha$
 in $\Phi^{i+1}_k$. The appropriate modification of the
decomposition given in Theorem \ref{seqlin} proves the claim in type $D_n$.
\end{proof}

\section{The Main Theorems}\label{main theorems section}

We now give a paving by affines for many Hessenberg varieties by
intersecting the Hessenberg variety with the cells of 
an appropriately chosen Schubert decomposition of $G/B$.
The first section describes Hessenberg varieties in
type $A_n$ and uses certain subsets of positive roots to parametrize
the cells of the paving and to give their dimension. 
The second section also discusses the paving of 
Hessenberg varieties in type
$A_n$ but describes this paving
using the combinatorics of Young tableaux rather than root systems.  
Section \ref{combi and Peterson} contains an extended description of the
Peterson variety to demonstrate how these results can be used computationally.
The third and fourth sections describe 
semisimple and regular nilpotent
Hessenberg varieties in classical types.

We make several comments which apply to all of the following results.
First, the choice of Schubert decomposition for $\H(M,H)$ is independent
of $M$.  In particular, the partial order on Hessenberg spaces gives rise
to 
filtrations of each Schubert cell into affine subspaces via 
$B \pi B \cap \H(M,H) \subseteq B \pi B \cap \H(M,H')$ if $H 
\subseteq H'$.

That each cell $B \pi B \cap \H(M,H)$ of the paving 
is an affine space over
the base field $\mathbb{C}$ implies the following result. 

\begin{proposition}
There is no odd-dimensional cohomology for Hessenberg varieties in type $A_n$
and for semisimple and regular nilpotent Hessenberg varieties in the other
classical types.
\end{proposition}

Note that while this paving can help identify the dimension of various
Hessenberg varieties, the question of whether all nilpotent 
Hessenberg varieties are equidimensional remains open.

\subsection{Hessenberg Varieties in Type $A_n$}\label{type A results}

We begin with nilpotent Hessenberg varieties and build to general
Hessenberg varieties in type $A_n$.  

\begin{theorem} \label{type A paving}
Let $G$ be $GL_n(\mathbb{C})$ or $SL_n(\mathbb{C})$
and let $N$ be a fixed nilpotent in $\n{}$.  Let $B'$ be the Borel
subgroup constructed by considering all upper triangular matrices in $G$ 
with respect to a basis which puts $N$ in Jordan canonical form.  

There
exists a permutation $\sigma$ such that the Borel $B = \sigma^{-1} \cdot B'$
induces a Bruhat 
decomposition whose Schubert cells intersect each Hessenberg variety
$\H(N,H)$ in a paving by affines.  The nonempty cells are $B \pi B$ with 
$\pi^{-1} \cdot N \in H$ and have dimension 
\[|\Phi_{\pi}| - |\C{\pi \cdot H} \cap \Phi_{U_{\pi} \cdot N}|.\]
\end{theorem}

\begin{proof}
To construct this Bruhat decomposition, first fix a basis with respect to
which $N$ is in Jordan canonical form.  Order the Jordan blocks from 
smallest to largest, fixing an order 
among equal-dimensional Jordan blocks once and for all.  
Then permute this basis according to the following
rules.   Index the basis vectors in $\ker N$ from $e_1$ to
$e_{|\ker N|}$ according to the order of the Jordan blocks containing $e_i$
so that $e_1$ belongs to the smallest Jordan block and $e_{|\ker N|}$ to
the largest Jordan block.  Given 
an ordering of the basis vectors in $\ker N^j$, index the
basis vectors in $\ker N^{j+1}/ \ker N^j$ from 
$e_{|\ker N^j|+1}$ to $e_{|\ker N^{j+1}|}$ increasing the index 
according to the
dimension of the Jordan block in which each $e_i$ is contained.

To show that $\Phi_N$ satisfies the desired conditions, we need to
describe more precisely this basis.  Suppose that the nilpotent $N$ 
lies in the conjugacy class
corresponding to the partition $\mu = (\mu_1, \ldots, \mu_s)$ 
with the parts ordered so that $\mu_i \geq \mu_{i+1}$.
Let $\mu' = (\mu_1', \ldots, \mu_{s'}')$ be the dual partition
associated to $\mu$.  

Define a function 
$w: \{1, \ldots, n\} \longrightarrow \{1, \ldots, s'\}$
by
\[ \sum_{k=1}^{w(i)-1} \mu_k' \leq i < \sum_{k=1}^{w(i)} \mu_k'  \hspace{1em}
        \textup{  for all  } i \in \{1, \ldots, n\} .\]
Note that $w$ is a nonincreasing function, that is $w(i) \geq w(i+1)$.

Define the roots
\[\beta_i = \sum_{j=0}^{\mu_{w(i)}' - 1} \alpha_{i-j}\]
for $\mu_1' \leq i < n$.  These roots have the following two properties:
\begin{enumerate}
\item Either $\beta_i - \alpha_i \in \Phi^+ \mbox{  or  } \beta_i = \alpha_i$.
\item $\beta_i \geq \alpha_j \Longrightarrow j \leq i$.
\end{enumerate}
In addition, the sequence $\{\h \beta_i\}$ is nonincreasing.
Consequently, the set $\{\beta_i: \mu_1' \leq i < n\}$ is 
non-overlapping.

There is 
a root space decomposition of $\g{}$ so that
$N = \sum_{i=\mu_1'}^{n-1} E_{\beta_i}$.
In fact, writing $\g{}$ with respect to the basis specified at the outset
of the proof gives such a root space decomposition.  By Lemma 
\ref{Hess invariant of orbit}, the isomorphism class of the 
Hessenberg variety is invariant under the choice of basis.  
The set $\Phi_N$ equals $\{\beta_i\}$ for this decomposition.

To see that $\Phi_N$ satisfies the conditions of Corollary 
\ref{constantrank on cell}, note that if $\h \beta_{i+k-1} < k$ then
$\Phi^i_{k,-N} = \Phi^i_j$ for $j = k - \h \beta_{i+k-1}$.  If not then
$\Phi^i_{k,-N}$ is empty since $\Phi_N$ is non-overlapping.  Moreover, 
for $\alpha$ to be in $\Phi_{U \cdot N} \cap \Phi^i_k$ means that
there exists $\beta_j$ in $\Phi_N$ with $\alpha \geq \beta_j$.
Since $\Phi_N$ is non-overlapping, we further conclude that
$\alpha \geq \beta_{i+k-1}$.  It follows that $\{\alpha\}^i_{-N}$ is 
nonempty if $\alpha$ is in $\Phi^i \cap \Phi_{U \cdot N} \cap \C{\pi \cdot H}$.
%%To see that $\Phi_N$ satisfies the conditions of corollary
%%\ref{constantrank on cell}, choose an element $u$ in $U$ and consider any root
%% $\beta_i + \gamma$ in
%%$\Phi_{u^{-1} \cdot N - N}$, for a nontrivial sum $\gamma$ 
%%of positive roots.
%%By hypothesis, $\beta_i + \gamma$ has an extremal simple root $\alpha_{i'}$ for
%%$i' \geq i$.  If $i=i'$, then $\beta_i$ is extremal for 
%%$\beta_i + \gamma$.  Furthermore, $\gamma$ is in the same row as
%%$\beta_i + \gamma$, since they agree in their smaller-indexed 
%%extremal simple roots.  If $i' > i$, then
%%\[\h(\beta_i + \gamma) > \h(\beta_i) \geq \h(\beta_{i'}).\]
%%So $\beta_{i'}$ is extremal for $\beta_i + \gamma$, and
%%$\beta_i + \gamma - \beta_{i'}$ is in the same row
%%as $\beta_i + \gamma$.
%%
%%This shows that $\Phi_N$ is saturated in $\Phi_{U \cdot N - N}$
%%and that 
%%\[\{ \alpha- f(\alpha): \alpha \in \Phi^i-\{\alpha_i\}\} \subseteq \Phi^i\]
%% for each $i$.
%%
%%Take $\lambda = (\Phi^n, \ldots, \Phi^1)$ and $F= (F_n, \ldots, F_1)$,
%%where
%%\[F_i = \{ \rho_{\alpha}(u^{-1} \cdot N): \alpha \in \Phi^i\}.\]

Either the hypotheses of Corollary \ref{empty cells}
hold or those of Lemma \ref{constantrank on cell} hold.
The claims then follow from Theorem \ref{affines}.
\end{proof}

This result extends to general linear operators $M$ by use of Theorem
\ref{reduce to nilp}.

\begin{theorem}
Given a linear operator $M$, there exists a Bruhat decomposition whose
intersection with each Hessenberg variety $\H(M,H)$
is a paving by affines of $\H(M,H)$.  Its
nonempty cells are $B \pi B \cap \H(M,H)$ 
such that $\pi^{-1} \cdot M \in H$ and have dimension
\[|\Phi_{\pi}| -|\Phi_{U_{\pi,M} \cdot (M-s)} 
                \cap \C{\pi \cdot H}| - 
        |\Phi_{\pi} \cap \Phi_{\n{M}} \cap \C{\pi \cdot H}|.\]
\end{theorem}

\begin{proof}
Choose a basis that puts $M$ in Jordan canonical form.  Permute the 
basis vectors within each generalized eigenspace so that the nilpotent
part of $M$ on each Jordan block is in the form specified by Theorem
\ref{type A paving}.  Write $M = \sum (S_j + N_j)$ where $S_j$ is
diagonal, $N_j$ is nilpotent, and $S_j + N_j$ is the $j^{th}$ Jordan 
block of $M$.  Again using Lemma 
\ref{Hess invariant of orbit}, we see that the 
isomorphism class of the Hessenberg variety $\H(M,H)$ is preserved
under this change of basis.

If $\Delta_j$ are the simple roots whose root vectors generate the $j^{th}$
Jordan block then $\Phi^+_{\sum S_j} = \sp \langle \bigcup \Delta_j \rangle$.
In particular, each $S_j + N_j$ is non-overlapping since $\Phi_{N_j}$ is 
contained in $\sp \Delta_j$.  
The conditions of Part (1) in Theorem
\ref{reduce to nilp} are thus met.

As usual, use $\mathfrak{p}_{\Delta_j}$ to denote the parabolic subalgebra
associated to $\Delta_j$ and 
$\mathfrak{l}_{\Delta_j}$ to 
denote its Levi subalgebra.  
Write $\mathfrak{p}_M$ for the 
parabolic associated to $\bigcup \Delta_j$ and $\n{M}$ (respectively
$\mathfrak{l}_M)$ 
for its nilpotent (respectively Levi) part.  All Levi
subalgebras are assumed to contain the diagonal matrices.

By Theorem \ref{type A paving}, the variety $U(N_j,\mathfrak{l}_{\Delta_j}
\cap \n{\pi \cdot H}) \cap U_{\pi} \cap U_{\Delta_j}$ 
is nonempty only if $\pi^{-1} \cdot N_j$
is in $H$.  In that case the variety is affine space of dimension
$|\Phi_{\pi} \cap \sp \Delta_j| - | \Phi_{U_{\Delta_j}\cap U_{\pi} 
\cdot N_j} \cap \C{\pi
\cdot H}|$.  This ensures the conditions of Part (2) of Theorem 
\ref{reduce to nilp}.  The claim follows by summing over $j$.  
\end{proof}

\subsection{Combinatorial Description of the Paving}\label{combi and Peterson}

Let $E_{ij}$ be the standard basis for $\mathfrak{gl}_n(\mathbb{C})$ 
in which $E_{i,i+1}$
corresponds to the root vector $E_{\alpha_i}$.  
This basis can be used to construct a bijection between the set of
Hessenberg spaces $H$
and  functions $h: \{1,\ldots, n\} \longrightarrow \{1,\ldots, n\}$
such that 
\[h(i) \geq \max \{i,h(i-1)\}.\]
Explicitly, the element $E_{ij}$ is in
$H$ if and only if $i \leq h(j)$.  We use this basis throughout
this section to establish a bijection between the dimension of the cells
from Section \ref{type A results} and
the number of certain configurations in specific Young tableaux.

Let $N$ be a nilpotent whose conjugacy class corresponds to the partition
$\mu = (\mu_1, \ldots, \mu_s)$ with the parts ordered so that
$\mu_i \geq \mu_{i+1}$.  Associate to $N$ the Young diagram whose $i^{th}$ 
column has $\mu_i$ blocks.  We use the convention that the blocks in this 
Young diagram are both left-aligned and bottom-aligned as in this 
example.
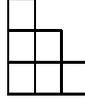
\[\begin{tabular}{|c|c|c|}
\cline{1-1}   & \multicolumn{2}{c}{  } \\
\cline{1-2}   &  & \multicolumn{1}{c}{  } \\
\cline{1-3}   &  &   \\
\hline
\end{tabular}
\]
We refer to this as the Young diagram associated to $N$.

\begin{theorem}\label{type A combi paving}
The nonempty cells of the paving of $\H(N,H)$ given in Theorem 
\ref{type A paving} correspond to those fillings of the Young
diagram associated to $N$ for which the configuration
\[
\begin{tabular}{|c|}
\cline{1-1}  $\pi^{-1} k$  \\
\cline{1-1}  $\pi^{-1} j$  \\
\hline
\end{tabular}\]
only occurs if $\pi^{-1} j \leq h(\pi^{-1} k)$.

Given a nonempty cell represented as a (filled) Young tableau, the
dimension of this cell is the sum of the following two quantities.
\begin{enumerate}
\item The number of configurations
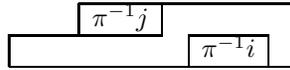
\[
\begin{tabular}{|cccccc|}
\cline{3-6}  \multicolumn{2}{c|}{  } & \multicolumn{1}{|c|}{$\pi^{-1} j$} & \multicolumn{2}{c}{   } & \\
\cline{1-3} \cline{5-5}  \multicolumn{2}{|c}{ $\hspace{1.6em}$ }& & &\multicolumn{1}{|c|}{$\pi^{-1} i$} &  \\
\hline
\end{tabular}\]
where box $i$ is to the right of or below box $j$, there is no box above $j$,
and the values filling these boxes satisfy $\pi^{-1} i > \pi^{-1} j$. 
\item The number of configurations
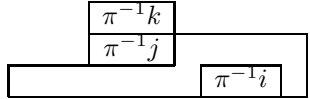
\[
\begin{tabular}{|cccccc|}
\cline{3-3}  \multicolumn{2}{c|}{$\hspace{2em}$  } & \multicolumn{1}{|c|}{$\pi^{-1} k$} & \multicolumn{3}{|c}{  } \\
\cline{3-6}  \multicolumn{2}{c|}{  } & \multicolumn{1}{|c|}{$\pi^{-1} j$} & \multicolumn{2}{c}{   } & \\
\cline{1-3} \cline{5-5}  \multicolumn{2}{|c}{ $\hspace{1.6em}$ }& & &\multicolumn{1}{|c|}{$\pi^{-1} i$} &  \\
\hline
\end{tabular}\]
where box $i$ is to the right of or below box $j$ and the values filling these
boxes satisfy $\pi^{-1} j < \pi^{-1} i \leq h(\pi^{-1} k)$.
\end{enumerate}
\end{theorem}

\begin{proof}
We begin by fixing a number in $\{1, \ldots, n\}$
to each box in the Young diagram associated to $N$.
The blocks are indexed from the bottom rightmost box to the 
top leftmost box by incrementing leftwards along each row then
going to the rightmost block on the next higher row and repeating as needed.
For instance, the Young diagram shown previously is indexed as follows:
\[
\begin{tabular}{|c|c|c|}
\cline{1-1}  6 & \multicolumn{2}{c}{  } \\
\cline{1-2}  5 & 4 & \multicolumn{1}{c}{  } \\
\cline{1-3}  3 & 2 & 1  \\
\hline
\end{tabular}\]
With respect to the standard matrix basis, the expression for $N$ given in 
Theorem \ref{type A paving} is equivalent to
\[ N = \sum_{\{(j,k): \mbox{ {\scriptsize box $k$ above box $j$}}\}} 
                E_{jk}.  \]
In other words, $N$ sends the basis vector $e_k$ to $e_j$ if and only if 
the $k^{th}$ box lies above the $j^{th}$ box.

Given this Young diagram, we can describe the Bruhat cells as
 (filled) Young tableaux.  In particular, 
we associate a Young tableau to each permutation $\pi$ by filling
the $i^{th}$ block with $\pi^{-1} i$.  The roots in $\Phi_{\pi}$ are
indexed by the set
\[ \Phi_{\pi} = \bigcup_{1 \leq j \leq n} \{(i,j): i < j, \pi^{-1} i > \pi^{-1} j\}.\]
This is the number of boxes to the right of or below the $j^{th}$ box that
are filled with numbers greater than that filling the $j^{th}$ box.

A basis vector $E_{ij}$ is in $H$ if and only if $h(j) \geq i$.  
Consequently, the set $\C{\pi \cdot H}$ can be characterized as
\[ \C{\pi \cdot H} = \{(\pi i,\pi k): i > h (k)\}. \]
Lemma \ref{sat gives rankcount} showed that the roots in 
$\Phi_{U_{\pi} \cdot N} \cap \C{\pi \cdot H}$ correspond bijectively
to the elements of the set
\[\{(i,k): i<j, \pi^{-1} i > \pi^{-1} j, \mbox{box $j$ is below box $k$} \}\cap \C{\pi \cdot H}.\]
This is because each root in $\Phi_{U_{\pi} \cdot N}$ can be written as
a sum $\alpha + \beta$ with $\beta$ in a higher row than $\alpha$ and
with $\beta$ in $\Phi_N$.
The cardinality of this set is the same as that of
\[\{(i,k): i<j, \pi^{-1} i > \pi^{-1} j, \mbox{box $j$ is below box $k$}, \pi^{-1} i > h(\pi^{-1} k) \}.\]
It follows that the quantity
$|\Phi_{\pi}| - |\Phi_{U_{\pi} \cdot N} \cap \C{\pi \cdot H}|$ is the 
cardinality of the set
\[\begin{array}{c} \{(i,k): i<j, \pi^{-1} i > \pi^{-1} j, \mbox{ box $j$ is below no box} \} 
\cup \\
\{(i,k): i<j, \pi^{-1} i > \pi^{-1} j, \mbox{box $j$ is below box $k$}, 
\pi^{-1} i \leq h(\pi^{-1} k) \}.\end{array}\]

Given a Young tableau, the cardinality of the first set is the same 
as the number of configurations
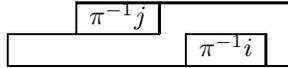
\[
\begin{tabular}{|cccccc|}
\cline{3-6}  \multicolumn{2}{c|}{  } & \multicolumn{1}{|c|}{$\pi^{-1} j$} & \multicolumn{2}{c}{  } & \\
\cline{1-3} \cline{5-5}  \multicolumn{2}{|c}{ $\hspace{1.6em}$ }& & &\multicolumn{1}{|c|}{$\pi^{-1} i$} &  \\
\hline
\end{tabular}\]
where box $i$ is to the right of or below box $j$, there is no box above $j$,
and the values filling these boxes satisfy $\pi^{-1} i > \pi^{-1} j$.  The cardinality
of the second set is the number of configurations
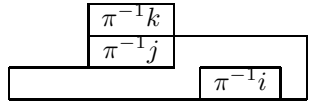
\[
\begin{tabular}{|cccccc|}
\cline{3-3}  \multicolumn{2}{c|}{$\hspace{2em}$  } & \multicolumn{1}{|c|}{$\pi^{-1} k$} & \multicolumn{3}{|c}{  } \\
\cline{3-6}  \multicolumn{2}{c|}{  } & \multicolumn{1}{|c|}{$\pi^{-1} j$} & \multicolumn{2}{c}{   } & \\
\cline{1-3} \cline{5-5}  \multicolumn{2}{|c}{ $\hspace{1.6em}$ }& & &\multicolumn{1}{|c|}{$\pi^{-1} i$} &  \\
\hline
\end{tabular}\]
where box $i$ is to the right of or below box $j$ and the values filling these
boxes satisfy $\pi^{-1} j < \pi^{-1} i \leq h(\pi^{-1} k)$.

Finally, the Schubert cell corresponding to $\pi$ is nonempty if and only if
$\pi^{-1} \cdot N$ is in $H$.  This is equivalent to the statement that
$\pi^{-1} j \leq h(\pi^{-1} k)$ for each instance when 
the $j^{th}$ box is below the $k^{th}$ box.
\end{proof}

The following result demonstrates how this theorem can be used computationally.
It gives the Betti numbers for the Peterson variety.  The Peterson variety
was defined generally in the Introduction.  In type $A_n$, it can be more
simply described.  Let $N$ be a regular nilpotent operator on $\mathbb{C}^n$,
that is a nilpotent with a single Jordan block.
Let $V = (V_1, \ldots, V_n)$ denote a full flag in $\mathbb{C}^n$, which is to
say that each $V_i$ is an $i$-dimensional complex vector space and $V_i
\subseteq V_{i+1}$ for $i=1, \ldots, n-1$.  The Peterson variety $\H(N,H)$
is defined to be
\[\H(N,H) = \{V: V \textup{  is a full flag, }NV_i \subseteq V_{i+1} 
    \textup{  for  } i=1, \ldots, n-1\}.\]
This corresponds to choosing 
\[H = \sp \langle E_{\alpha}: -\alpha \in \Phi^-
\cup \Delta \rangle = \sp \langle E_{ij}: 1\leq i,j \leq n, i\leq j+1 \rangle\]
 or, equivalently, choosing $H$ to be the set of $n \times
n$ matrices which are zero below the subdiagonal.

If the basis for $\mathbb{C}^n$ is chosen so that $N$ is in Jordan canonical
form, the flags in the Peterson variety can be described completely as
follows.  First consider matrices of the form
\[\left( \begin{array}{cccc} c & b & a & 1 \\ b & a & 1 & 0 \\
a & 1 & 0 & 0 \\ 1 & 0 & 0 & 0 \end{array} \right)\]
with ones along the cross-diagonal, zeroes below, and constant values
along each line above and parallel to the cross-diagonal.  Write $J_i$ for
an $i \times i$ matrix of this form.  The flags in $\H(N,H)$ are in bijective
correspondance to matrices of the form 
\[\left( \begin{array}{cccc}J_{i_1}
& 0 & 0 & 0\\ 0 & J_{i_2} & 0 & 0 \\ 0 & 0 & \ddots & 0 \\ 0 & 0 & 0 & J_{i_k}
\end{array} \right).  \]
To obtain a flag from a matrix, let $V_i$ be the
span of the first $i$ columns.

\begin{theorem}
Let $N$ be a regular nilpotent.  Let $H$ be the subspace of $\g{}$ given by
\[H = \mathfrak{b} \oplus \bigoplus_{\alpha \in \Delta} \g{-\alpha}.\]

There is a natural bijection between the cells of the Peterson variety 
$\H(N,H)$ and the ordered
partitions of $n$.  If $C = (c_1, \ldots, c_k)$ is the cell whose associated
partition has parts $c_i$ then the complex dimension of $C$ is
$n-k$. The number of cells with complex dimension $k$ is 
$\binom{n-1}{k}$.
\end{theorem}

\begin{proof}
A regular nilpotent has only one Jordan block so the associated Young diagram
consists of a single column with $n$ boxes.  
Since $h(i) = i+1$, Theorem
\ref{type A combi paving} tells us that
 each adjacent pair of boxes 
\[
\begin{tabular}{|c|}
\cline{1-1} $\pi^{-1} (i+1)$  \\
\cline{1-1} $\pi^{-1} (i)$  \\
\hline \end{tabular} \]
in a Young tableau for a nonempty cell 
satisfies either $\pi^{-1} (i+1) = \pi^{-1} (i)-1$ or $\pi^{-1} (i+1) > \pi^{-1} (i)$.

In fact, we can characterize the Young tableaux completely by selecting
boxes of the column to initiate increasing runs 
\[(\pi^{-1}(i),\pi^{-1}(i)+1, \ldots, \pi^{-1}(i)+c-1)\]
 of the numbers
filling the boxes.  Suppose that $(i_1, \ldots, i_1+c_1-1)$
and $(i_2, \ldots, i_2+c_2-1)$ index the boxes in two increasing runs.
If the Young tableau represents a nonempty cell, then $\pi^{-1}(i_1)<\pi^{-1}(i_2)$ implies 
that the first run must fill lower boxes of the Young diagram than the second. 
The sizes of the increasing runs give an ordered partition of $n$.

In addition, given an arrangement like
\[
\begin{tabular}{|c|}
\cline{1-1} $\pi^{-1} (j+1)$  \\
\cline{1-1} $\pi^{-1} (j)$  \\
\cline{1-1} $\vdots$ \\
\cline{1-1} $\pi^{-1} (i)$ \\
\hline \end{tabular} \]
the quantity $\pi^{-1}(i)$ is greater than $\pi^{-1}(j)$ only if $\pi^{-1}(j)$ is in the same
increasing run as $\pi^{-1}(i)$.  If $\pi^{-1}(j+1)$ is also in this increasing run
then $h(\pi^{-1}(j+1)) = \pi^{-1}(j+1)+1 = \pi^{-1}(j)$ so the arrangement does not
contribute to the dimension of the cell, as per 
Theorem \ref{type A combi paving}.  
By comparing each lower box in the run to the topmost in that run, we see that
each
increasing run of length $c$ contributes exactly $c-1$ to the total dimension
of the cell.

Thus, if $C = (c_1, \ldots, c_k)$ is the cell whose associated
partition has parts $c_i$ then the dimension of $C$ is
\[ \dim C = \sum_{i=1}^k (c_i - 1) = n-k.\] 
The number of cells with $n-k$ parts is the same as the number of the
ways to choose $n-k-1$ out of the $n-1$ lower boxes of the Young diagram,
since the top box always initiates an increasing run.
\end{proof}

In \cite{ST} we give further analyses of the Poincare polynomials of
regular nilpotent Hessenberg varieties.

Semisimple and general Hessenberg varieties can also be
described combinatorially.  Associate to a linear operator $M$
the multidiagram with one Young diagram for each generalized
eigenspace of $M$.  The operator $M$ acts on each generalized eigenspace
by the sum of a nilpotent operator $N_j$ and a semisimple 
operator constant on the generalized eigenspace.  
The Young diagram corresponding to the $j^{th}$ 
generalized eigenspace is simply that associated to $N_j$.

Order the Young diagrams from largest to smallest,
right to left.  Index the boxes of each
Young diagrams as if for $N_j$.  Finally,
increment the indices of the $j^{th}$ Young diagram by
the number of boxes in the Young diagrams to the right of it.
For example, 

\begin{center}
\begin{tabular}{|c|c|}
\cline{1-2} 7 & 6 \\
\hline
\end{tabular}
\hspace{2em}
\begin{tabular}{|c|}
\cline{1-1} 5 \\
\cline{1-1} 4 \\
\hline
\end{tabular}
\hspace{2em}
\begin{tabular}{|c|c|}
\cline{1-1} 3 & \multicolumn{1}{c}{} \\
\cline{1-2} 2 & 1 \\
\hline
\end{tabular}
\end{center}

As before, describe the permutation $\pi$ 
by filling box $i$ with the value $\pi^{-1} i$.

\begin{theorem}
Given $M$ and its associated multidiagram ${\mathcal Y}$, 
there is a bijective correspondance 
between nonempty affine cells of $\H(M,H)$ and fillings
of ${\mathcal Y}$ such that the configuration
\begin{tabular}{|c|}
\cline{1-1} $\pi^{-1} i$ \\
\cline{1-1} $\pi^{-1} j$ \\
\hline \end{tabular}
occurs only if $\pi^{-1} j \leq h(\pi^{-1} i)$. 
The dimension of the affine corresponding to a permissable
filling of ${\mathcal Y}$ is the sum of the following 
quantities.
\begin{enumerate}
\item The number of configurations of type
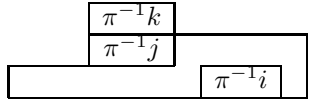
\[
\begin{tabular}{|cccccc|}
\cline{3-3}  \multicolumn{2}{c|}{$\hspace{2em}$  } & \multicolumn{1}{|c|}{$\pi^{-1} k$} & \multicolumn{3}{|c}{  } \\
\cline{3-6}  \multicolumn{2}{c|}{  } & \multicolumn{1}{|c|}{$\pi^{-1} j$} & \multicolumn{2}{c}{   } & \\
\cline{1-3} \cline{5-5}  \multicolumn{2}{|c}{ $\hspace{1.6em}$ }& & &\multicolumn{1}{|c|}{$\pi^{-1} i$} &  \\
\hline
\end{tabular}\]
where $i$ is less than $j$, box $i$ and box $j$ are in the
same Young diagram, and the values filling these
boxes satisfy $\pi^{-1} j < \pi^{-1} i \leq h(\pi^{-1} k)$.  
(If box $k$ does not exist, the latter inequality is 
considered vacuously satisfied.)
\item The number of 
configurations of type
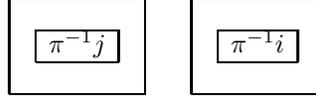
\begin{tabular}{|c|c|c|}
\cline{1-3} \multicolumn{3}{|c|}{$\hspace{2em}$} \\
\cline{2-2} & $\pi^{-1} j$ & \\
\cline{2-2} \multicolumn{3}{|c|}{$\hspace{2em}$} \\
\hline \end{tabular} $\hspace{1em}$
\begin{tabular}{|c|c|c|}
\cline{1-3} \multicolumn{3}{|c|}{$\hspace{2em}$} \\
\cline{2-2} & $\pi^{-1} i$ & \\
\cline{2-2} \multicolumn{3}{|c|}{$\hspace{2em}$} \\
\hline \end{tabular}
\newline with $h(\pi^{-1} j) \geq \pi^{-1} i > \pi^{-1} j$ and boxes $i$ and
$j$ in different Young diagrams.
\end{enumerate}
\end{theorem}

\begin{proof}
The root $\alpha = \alpha_i + \ldots + \alpha_{i+j}$ is in 
$\Phi_{\mathfrak{l}_{M}}^+$ if and only if the $i^{th}$ and 
$(i+j+1)^{th}$ boxes are in the same Young diagram.  This
root contributes one affine dimension to the total dimension
of the cell if and only if the conditions of Theorem \ref{type A combi paving}
are satisfied.

By contrast, $\alpha$ is in $\Phi_{\n{M}}$ exactly when the 
$i^{th}$ and $(i+j+1)^{th}$ boxes are in different Young 
diagrams.  By Theorem \ref{reduce to nilp}, this 
root contributes one dimension to the
total dimension of the cell only if $\alpha$ is in 
$\Phi_{\pi} \cap \pi \M_{H}$, that is if 
\[\pi^{-1} i > \pi^{-1} j \mbox{   and   } h(\pi^{-1} j) \geq \pi^{-1} i.\]
\end{proof}

\subsection{Semisimple Hessenberg Varieties in Classical Type}

Let $S$ be any se\-mi\-simple element in $\g{}$, a Lie algebra of
classical type.  Then the corresponding Hessenberg variety
$\H(S,H)$ can be paved by affines.  Once a root decomposition
is fixed, we establish the notation that $\mathfrak{p}$ denotes
the parabolic subalgebra generated by $\Phi^+_S$ and that $\n{}$ 
(respectively
$\mathfrak{l}$) denotes its nilpotent (respectively Levi) part.
We assume that $\mathfrak{l} \supseteq \mathfrak{h}$.

\begin{theorem}
If $S$ is a semisimple element in $\g{}$ a Lie algebra of
classical type then there exists a Bruhat decomposition
whose intersection with each Hessenberg
variety $\H(S,H)$ is a paving by affines of $\H(S,H)$.  
Each cell $B \pi B \cap \H(S,H)$ is nonempty and
has dimension
\[|\Phi_{\pi} \cap \pi \M_H \cap \Phi_{\n{}}| + 
      |\Phi_{\pi} \cap \Phi_{\mathfrak{l}}|.\]
\end{theorem}

\begin{proof}
Embed the corresponding group $G$ into $GL_N(\mathbb{C})$
in the natural way so that $\rk G = \lfloor N/2 \rfloor$
and so that the simple roots $\alpha_i$ for $G$ correspond
to those for $GL_N(\mathbb{C})$ when 
$i < \lfloor N/2 \rfloor$.

There exists a basis for $\mathbb{C}^N$ with respect to which
$S$ is diagonal and in Jordan canonical form in 
$\mathfrak{gl}_N(\mathbb{C})$.  In other words, the diagonal 
entries of $S$ are grouped by eigenvalue.  (If zero is an 
eigenvalue of $S$, we further assume that the diagonal of $S$
is zero in row $\lfloor N/2 \rfloor$.)  Write $S_0$ for this
element.  

Since $G \cdot S = GL_N(\mathbb{C}) \cdot S \cap \g{}$,
there exists an element $g$ in $G$ such that $g \cdot S = S_0$
(see \cite[section IV.2.19]{SS}).
By Lemma 
\ref{Hess invariant of orbit}, the 
Hessenberg variety $\H(S,H)$ is isomorphic to
$\H(S_0, g \cdot H)$.
Note that $\{\alpha_i : \alpha_i(S_0) = 0\}$ 
generate a parabolic subalgebra $\mathfrak{p}$.  The
operator $\ad S_0$ is zero on the Levi subalgebra $\mathfrak{l}
\subseteq \mathfrak{p}$ and is nonzero on each $E_{\alpha}$ in 
the nilpotent subalgebra $\n{} \subseteq \mathfrak{p}$.  Thus,
this decomposition satisfies Theorem \ref{reduce to nilp}.  The
conclusions follow.
\end{proof}

\subsection{Hessenberg Va\-rie\-ties of Regular Nilpotent Elements 
in Classical Type}

One class of nilpotent Hessenberg varieties can be paved by affines
using the same methods.  For this result, $G$ is
 a linear algebraic group of classical type.

\begin{theorem}\label{reg nilps theorem}
Fix a regular nilpotent element $N$, let $\mathfrak{b}$ be the unique
Borel subalgebra with $N \in \mathfrak{b}$, and let $B$ be the Borel
subgroup corresponding to $\mathfrak{b}$.  The intersection of the Bruhat 
decomposition with respect to $B$ and the Hessenberg variety $\H(N,H)$ 
is a paving by affines of $\H(N,H)$ for each $H$.  
The nonempty cells of this paving are $B \pi B \cap \H(N,H)$ satisfying 
$\pi^{-1} \cdot N \in H$ and have dimension 
\[|\Phi_{\pi}| - |\C{\pi \cdot H} \cap \Phi_{U_{\pi} \cdot N}|.\]
\end{theorem}

\begin{proof}
Given any regular nilpotent $N$ of classical type, there
exists a root space decomposition for $\g{}$ such that
\[N = \sum_{\alpha \in \Delta} E_{\alpha}\]
by \cite[sections 5.2 and 5.4]{CM}.  The Borel subalgebra $\mathfrak{b}$
corresponding to the positive roots in this decomposition is the unique Borel
subalgebra containing $N$ as shown in \cite[3.2.13 and 3.2.14]{CG}. By
Lemma \ref{Hess invariant of orbit}, 
the choice of root space decomposition preserves the
isomorphism class of the Hessenberg variety.
Note that $\Phi_N = \Delta$ is trivially non-overlapping
and that the rows $\Phi^i$ in classical type are vertical.  

Since $\Phi_N = \Delta$, each $\Phi^i_{k,N}$ is 
$\Phi^i_{k-1}$ unless $k=1$.  Only in type $D_n$ is 
there a set $\Phi^i_k$ of size two.  In this case, $\Phi^i_k = 
\Phi^i_{k+1,-N} = \Phi^i_{k-1,+N}$.  Furthermore, if
$\alpha = \sum_{\beta_j \in \Phi^i} \beta_j$ is in $\Phi^i_k$
then $\{\alpha\}_{-N} = \Phi^i_{k-1}$ which by row-verticality
 must be at least as large as each $\beta_j$.

By Corollary \ref{empty cells}, if $N$ is not in $\pi \cdot H$
then $U(N,\n{\pi \cdot H})$ is empty.  
Otherwise, all the hypotheses of Lemma 
\ref{constantrank on cell} are satisfied.  Finally, by Theorem
\ref{affines}, the dimension of the cell $\H(N,H) \cap B \pi B$ 
is precisely that of $U_{\pi} \cap
U(N,\n{\pi \cdot H})$.
\end{proof}

For a discussion of related results for regular nilpotent Hessenberg
varieties including a simpler dimension formula, the Euler characteristic,
and some properties of the Poincare polynomials, see \cite{ST}.

%%%\newpage
%%%
%%%\addcontentsline{toc}{chapter}{Appendix: List of Symbols}
\section{Appendix: List of Symbols}
%%%{\bf {\huge Appendix: List of Symbols}}

\begin{center}
\begin{tabular}{llll}
$G \hspace{.25in}$ & \pageref{G}$\hspace{2in}$ & 
$\Phi^i_k \hspace{.25in}$ & \pageref{row height partition} $\hspace{1in}$\\ 
$\g{}$ & \pageref{lie algebra} & 
$U_i$ & \pageref{group row} \\
$\mathfrak{h}$ & \pageref{cartan subalgebra} &
$U_{\pi}$ & \pageref{unipotent schubert cell} \\ 
$\Phi^+$ &   \pageref{positive roots} & 
$\rho_K$ & \pageref{projection to K} \\
$\Delta$ & \pageref{simple roots} & 
$\rho_i$ & \pageref{projection to row i} \\ 
$\g{\alpha}$ & \pageref{root space} &
$\rho_{\alpha}$ & \pageref{projection to alpha} \\
$E_{\alpha}$ & \pageref{root vector} & 
$\iota_K$ & \pageref{inclusion from K} \\
$\mathfrak{b}$ & \pageref{borel subalgebra} &
$\gamma_i$ & \pageref{long root row i} \\
$\n{}$ & \pageref{nilradical} &
$\Phi_M$ & \pageref{roots of M} \\
$H$ & \pageref{H} & 
$\Phi_{\mathfrak{u}}$ & \pageref{roots of mathfrak u}\\ 
$\M_H$ & \pageref{Hessenberg roots} & 
$\Phi_{\pi}$ & \pageref{roots of schubert cell} \\
$g \cdot M$ & \pageref{g cdot M} &
$\C{K}$ & \pageref{cal C: complement} \\
$\geq$ & \pageref{root order} & 
$U(M,\n{K})$ & \pageref{unipotent seq lin} \\
$\Phi^i$  & \pageref{row roots}&  
$\Phi_S^+$ & \pageref{roots annihilating S} \\ 
$\n{i}$ & \pageref{lie algebra row} &
$P^i_{\pm N}$ & \pageref{row i pm roots of N} \\
\end{tabular}
\end{center}

\end{document}